\documentclass[11pt]{article}
\pdfoutput=1
\usepackage[margin=1in]{geometry}          
\usepackage{graphicx}
\usepackage{amsthm, amsmath, amssymb}
\usepackage{setspace}
\usepackage{hyperref}
\usepackage{booktabs} 
\usepackage{multicol}
\usepackage{multirow}
\usepackage{subcaption}
\usepackage{epstopdf}
\usepackage{float}
\graphicspath{{figdata/}}

\usepackage[loose,nice]{units} 

\author{T. A. J. Ouermi \\
        University of Utah School of Computing \\
        U of U Scientific Computing and Imaging Institute \\
        \texttt{touermi@cs.utah.edu}\\
        \and
        Robert M. Kirby \\
        University of Utah School of Computing \\
        U of U Scientific Computing and Imaging Institute \\
        \texttt{kiby@cs.utah.edu}\\
        \and
        Martin Berzins \\
        University of Utah School of Computing \\
        U of U Scientific Computing and Imaging Institute \\
        \texttt{mb@sci.utah.edu}\\
        }
\date{Revised version June 15, 2022}

\title{Numerical Testing of a New Positivity-Preserving Interpolation Algorithm}

\begin{document}
\maketitle

\begin{abstract}
An important component of a number of computational modeling algorithms is an interpolation method that preserves the positivity of the function being interpolated. 
This report describes the numerical testing of a new positivity-preserving algorithm that is designed to be used when interpolating from a solution defined on one grid to different spatial grid.  
The motivating application for this work was 
a numerical weather prediction (NWP) code that uses a spectral element mesh discretization for its dynamics core
and a cartesian tensor product mesh for the evaluation of its physics routines.  
This coupling of spectral element mesh, which uses
nonuniformly spaced quadrature/collocation points, and uniformly-spaced cartesian mesh combined with the desire to maintain 
positivity when moving between these meshes necessitates our work. 
This new approach is evaluated against several typical algorithms in use on a range of test problems in one or more space dimensions.
The results obtained show that the new method is competitive in terms of observed accuracy while at the same time preserving the underlying positivity
of the functions being interpolated. 
\end{abstract}

\section*{Introduction}
\label{sec:intro}
Interpolating from one grid to another is a fundamental part of a number of computational problems.
Furthermore, when interpolating solution values, it is important in some applications to conserve properties such as non-negativity of the solution.
For example in weather forecasting,
when mapping between the different grids used to calculate the dynamics and those used to calculate the physics, the polynomial approximation of positive quantities such as mass, density, or cloud mixing ratio may introduce negative values that are nonphysical.
These negative values may cause incorrect representations of other calculations that are dependent on these approximations.
In other applications such as parallel resilience \cite{Damodar},  combustion simulations and the solution of hyperbolic equations with ENO and WENO schemes, 
it is important to have polynomial approximations that preserve positivity on general meshes \cite{Berzins2010c,Berzins2013b}. 
A somewhat stronger condition is to use methods that preserve the local boundedness of the solution \cite{Berzins} when interpolating.

The particular example motivating this experimental study is the physics-dynamics coupling (PDC) module in the Navy Environmental Prediction System
using the NUMA Core (NEPTUNE) \cite{NUMA,Neptune:Alex}.
NEPTUNE is a next-generation global NWP system being developed at the Naval Research Laboratory (NRL) and the Naval Postgraduate School (NPS) \cite{Neptune:Alex}.
In NEPTUNE \cite{NUMA,Neptune:Alex}, the dynamics calculations are done on a spectral element mesh, whereas the physics routines require values on a uniform mesh.
In the context of NEPTUNE and similar codes, interpolating the solution values produced by the dynamics routines to the spatial 
points used by the physics routines and vice versa may lead to negative values unless care is taken.   
This result was shown by Skamarock et al. \cite{skamrock} who demonstrated that not preserving positivity may 
lead to unphysical results in a predicted physical quantity of interest, such as moisture.
Moreover, the nonphysical values introduced through interpolation may 
lead to spurious values in the results, 
compared to those produced with positivity-preserving interpolation \cite{tajo2022PPIsoftware}.

This report is concerned with the numerical testing of a new interpolation algorithm that has been proposed for positivity preservation when mapping solution values between structured meshes. 
The theoretical basis for the algorithm builds on the data-bounded work of Berzins \cite{Berzins} to develop a new data-bounded and positivity-preserving methods for both evenly- and unevenly-spaced structure meshes.
The new data-bounded interpolation (DBI) method in ~\cite{ouermi2022eno} relaxes conditions for data-boundedness, which gives greater accuracy than the conditions used in ~\cite{Berzins}. 
Ouermi et al. ~\cite{ouermi2022eno} further extended the DBI method to give a new positivity-preserving interpolation (PPI) method.
The application of these new methods to numerical weather prediction examples is described in \cite{tajo2022PPIsoftware}.
In this report, a number of possible alternative interpolation schemes are introduced.
A representative sample of such methods is compared against the new approaches on 
a number of different test functions, including smooth, $C^{0}$, discontinuous, and steep-gradient functions. 
The comparison undertaken focuses on how accurately the different methods are able to represent this underlying set of test functions.
In addition, a representative weather model problem is considered. 
Overall, it will be shown that the new methods are well suited for function approximation and mapping data values between meshes for numerical weather examples.
The generality of this approach suggests that these methods also have application to other problems for which preserving positivity is important.
\section{Examples of Existing Interpolation Methods}
\label{sec:background}
This section highlights several approaches that have been developed to address the need for data-bounded,  positivity-preserving, and shape-preserving interpolation. 
While this selection of methods is not all-inclusive, it is intended to illustrate the main types of polynomial-based approaches.

\subsection{Cubic Splines}
In Computer-Aided Design (CAD), graphics, 
and visualization, significant contributions have been made to develop and advance shape-preserving methods.  
Many of the approaches for shape-preservation are based on cubic splines.
In \cite{Schmidt1987} and \cite{Schmidt1988}, Schmidt and He{\ss} introduced positive interpolation methods using rational quadratic and cubic splines 
respectively.
Necessary and sufficient conditions for positivity are provided for both the rational quadratic and cubic interpolants.
These conditions impose some restrictions on the values of the first derivatives at each node.
As in \cite{HE199451}, both approaches lead to multiple solutions, and the one with the minimal curvature is selected.
The work in \cite{abdul2016shape}, \cite{karim2015positivity},  and \cite{HUSSAIN2008446} 
presented positivity-preserving interpolation methods that 
rely on rational cubic splines.
The $C^{2}$ continuity in \cite{abdul2016shape} is obtained by solving a tridiagonal system of linear equations.
All three methods introduce free parameters that are used to derive and enforce conditions for positivity. 
Butt and Brodlie \cite{BUTT199355} provide a method for constructing $C^{1}$ cubic Hermite splines.
This method is dependent on the availability of values of first derivatives at the nodes,
which may not be available in practice.
Positivity is enforced by imposing a bound on the values of the derivatives.
In the case where bounds on the derivatives are not met, one or two knots are inserted to ensure that the constructed spline is positive.          
Perhaps the most widely used approach for preserving monotonicity in many applications is PCHIP 
by Fritch and Carlson \cite{fritsch1980monotone} who derived necessary and sufficient conditions for monotone cubic interpolation, 
and provided an algorithm for building a piecewise cubic approximation from data.
This algorithm calculates the values of the first derivatives at the nodes based on the necessary and sufficient conditions.

\subsection{Quartic and Quintic Splines}
Although many shape-preserving interpolation methods are  cubic or lower order, a number of approaches target 
higher-order interpolants, 
with an emphasis on quartic or quintic polynomial approximations.
The work in \cite{10.1093/imanum/23.2.175} and \cite{10.1093/imanum/drt072} presents geometric or visual continuity $G^{1}$ and $G^{2}$ 
continuous shape-preserving interpolation using Pythagorean-Hodograph quintic splines curves. 
This approach uses Bernstein basis functions and a parametric representation of the interpolant in each interval. 
A sufficient condition for shape preservation is constructed based on free angular parameters that influence the shape of the curve in each interval. 
The appropriate angular parameters are selected based on the cubic-cubic (CC) criterion introduced in \cite{FAROUKI2008274}.
The $G^{2}$ case requires a tridiagonal solve and use of a Newton-Raphson iteration, which potentially affects the computational performance.
%

Hussain et al. \cite{Hussain2018ACR} and Hussain et al. \cite{hussain2009c2} introduced $C^{2}$ rational quintic interpolation interpolation 
approaches that preserve positivity.
These rational quintic functions are constructed with free parameters that are used to enforce positivity.
Both methods require the approximation of values of first and second derivatives at the nodes if these derivatives are not available.
In addition, the rational quintic interpolation methods in \cite{Hussain2018ACR} and \cite{hussain2009c2} have a $O(h^3)$ order of accuracy.
%

He{\ss} and Schmidt \cite{HE199451} developed interpolation schemes that preserve positivity and monotonicity using $C^2$ quartic and quintic splines.
Positivity and monotonicity are achieved by imposing some restrictions on the values of the first and second derivatives at each node.
This approach leads to a potentially infinite number of solutions that meet the required conditions. 
Of these solutions,  the solution with minimal curvature is selected using global minimization. 
The global nature of the minimization makes the algorithm challenging to parallelize and may have an impact on computational performance.
%
%
%
MQS \cite{Lux2019ANAF} is an example of a monotonic quintic spline method that was developed by 
Lux et al. \cite{Lux2019ANAF} who built on the work of He{\ss} and Schmidt \cite{10.1007/BF01934097}, 
and Ulrich and Watson \cite{doi:10.1137/0915035}.
This algorithm uses the sufficient conditions from \cite{10.1007/BF01934097} to check for monotonicity and the work in 
\cite{doi:10.1137/0915035} to adjust values of the first and second derivatives to ensure monotonicity.
This method requires the values of the first and second derivatives at the nodes, which may not be available in practice.
In this report, the first and second derivatives are approximated using a fourth-order finite difference stencil 
based on \cite{fornberg1988generation}.

\subsection{SPS and B-spline Higher Order Splines}
Costantini \cite{10.1145/264029.264050, 10.1145/264029.264059} developed a $C^{1}$ and $C^{2}$ 
Shape-Preserving Spline (SPS)  interpolation method using Berstein-Bezier polynomials of an arbitrary degree.
The desired shape property is obtained by imposing restrictions on the value of the first derivatives at the nodes.
The Bezier coefficients for each spline are derived from a linear function.   
For a given interval, the coefficients of the Berstein-Bezier polynomial interpolant are selected from a linear function.
The first derivatives at the nodes are calculated such that the sufficient conditions for shape preservation given in \cite{10.1145/264029.264050} are met.
The approximation of the first derivatives at the node is third-order accurate.
In addition, Theorem 9 of \cite{Costantini1990} shows that the spline method presented in \cite{10.1145/264029.264050, 10.1145/264029.264059} has an error of $O(h^{4})$.
More details on the construction of the splines, an algorithm and a software package for the SPS method can be found in \cite{10.1145/264029.264050, 10.1145/264029.264059}.
%
In addition to the positivity-preserving approaches, conventional B-splines \cite{de1978practical} are also used here.
Although the B-spline approach does not preserve-positivity, many of the approaches mentioned in this work are based on B-splines
and so the use of unmodified B-splines provides an accuracy check on the other spline methods.

\subsection{DBI and PPI Methods}
\label{subsec:dbippi}
The numerical solution of partial differential equations (PDEs), particularly hyperbolic equations,
is an another area in which various methods have been developed to enable data-bounded and positivity-preserving approximations.
In order to preserve positivity in discontinuous Galerkin (dG) schemes, Zhang et al. 
\cite{ZHANG2017, Zhang2012, Zhang2752} and Light et al. \cite{light} introduced a 
linear rescaling of polynomials that ensures that the evaluation of the polynomial 
at the quadrature points is positive.  
In addition, this linear rescaling of the polynomial conserves mass. 
The polynomial rescaling, however, does not address the case of interpolating between different meshes, which is the primary focus of this work.
Harten et al. \cite{HARTEN19973} developed an essentially non-oscillatory (ENO) piece-wise polynomial reconstruction that is suitable 
for interpolating between different meshes.
ENO methods adaptively build an interpolant based on Newton divided differences and can help remove Gibbs-like effects but do not guarantee positivity.
A weighted combination of ENO schemes, (WENO) has been used by Zhang et al. \cite{Zhang2012_2} and many
others. 

A DBI method was developed by Berzins using evenly spaced meshes from ENO methods \cite{Berzins}.
This method was extended by the authors in ~\cite{ouermi2022eno,tajo2022PPIsoftware} to work for both evenly and unevenly spaced meshes and, more importantly, to the PPI method. 
Ouermi et al. ~\cite{ouermi2022eno} relaxed the conditions for data-boundedness which, gives greater accuracy compared to the conditions used in ~\cite{Berzins}.
Both the improved DBI and the new PPI methods are used in this report.
The PPI method further extends the DBI method by relaxing the bounds on the ratio of divided differences 
and so allows the interpolant to grow beyond the data, while still remaining positive.
For a given interval, the DBI and PPI methods successively select stencil points 
until the required bounds are violated or $d+1$ points are selected, with $d$ being the target degree 
of the interpolant.
In addition to enforcing data-boundedness and positivity, the algorithm in ~\cite{tajo2022PPIsoftware} uses a user-supplied parameter $st$ to guide the stencil construction procedure. 
When adding the next point to both the right or left of the current stencil meets the requirements for data-boundedness or positivity, the algorithm makes the selection based on the three cases below.  
\begin{itemize}
  \item If $st=1$, the algorithm chooses the point with the smallest divided difference, as in the ENO stencil. 
  \item If $st=2$, the point to the left of the current stencil is selected if the number of points to the left of $x_{i}$ is smaller than the number of points to right. 
  Similarly, the point to the right is selected if the number of points to the right of $x_{i}$ is smaller than the number of points to the left.
  When both the number of points to right and left are the same, the algorithm chooses the point with the smallest ratio of divided differences.
  \item If $st=3$, the algorithm chooses the point that is closest to the starting interval $I_{i}$.
\end{itemize}
Enforcing positivity alone may still lead to undesirables oscillations.
To address this limitation the algorithm, provides the parameters $\epsilon_{0}$ and $\epsilon_{1}$ that are used to impose an upper and lower bound for each interpolant.
For each interval $I_{i}$, the bounds are constructed using the parameters $\epsilon_{0}$ and $\epsilon_{1}$, and the data values $u_{i}$ and $u_{i+1}$.
Both the DBI and PPI methods and the algorithm details are described in ~\cite{tajo2022PPIsoftware} with numerical examples pertaining to NWP. 
\section{Comparison Methodology}
\label{sec:method}

\subsection{Compared Methods}
\label{subsec:comp_method}
The numerical experiments in this report use the PCHIP \cite{fritsch1980monotone}, MQS \cite{Lux2019ANAF}, SPS ~\cite{10.1145/264029.264050, 10.1145/264029.264059}, 
B-splines \cite{de1978practical}, the improved DBI \cite{Berzins}, and the new PPI methods \cite{ouermi2022eno}.
These methods are available as follows:

{\bf PCHIP:}
The version of the PCHIP algorithm used in this report is implemented in Fortran 90 and can be found at \url{https://people.sc.fsu.edu/~jburkardt/f_src/pchip/pchip.html}.

{\bf  MQS:} The method of Lux et al. \cite{Lux2019ANAF} is an example of a method for monotonic quintic splines.
The algorithm is implemented in Python3 and can be found \url{https://github.com/tchlux/papers/tree/master/\%5B2019-11\%5D_HPC_(quintic_spline)}.  

{\bf SPS:}  Costantini \cite{10.1145/264029.264050, 10.1145/264029.264059} introduced a high-order Shape-Preserving (monotonicity-, and convexity-preserving) Spline (SPS) method using Berstein-Bezier polynomials of arbitrary degree. 
The SPS method is implemented in the BVSPIS software package in Fortran 77 and is available from ACM as Algorithm 770 \cite{10.1145/264029.264059}
\url{https://dl.acm.org/action/downloadSupplement?doi=10.1145\%2F264029.264059&file=770.gz&download=true}.

{\bf B-splines:}
PPPACK, a Fortran 90 library that evaluates piecewise polynomial functions, including cubic splines. The original FORTRAN77 library is by Carl de Boor
\cite {de1978practical}. 
The package is available from \url{https://people.sc.fsu.edu/~jburkardt/f_src/pppack/pppack.html}.  

{\bf HPPIS:}  The DBI and PPI methods have been developed based on the theory and \\ algorithm in ~\cite{ouermi2022eno,tajo2022PPIsoftware}.
The software and implementation details can be found in ~\cite{tajo2022PPIsoftware}.

\subsection{Comparison Criteria}
\label{subsec:com_criteria}
The three steps outlined below are used to compare the different methods when used to approximate smooth and nonsmooth functions.
The errors are measured in a discrete approximation to the  $L^{2}$-error norm. 
  \begin{itemize}
    \item The first step consists of demonstrating that the various schemes preserve positivity for each of the test functions used.
      In addition, this step is used to show that a standard polynomial interpolation method does not guarantee positivity. 
    \item The second step experimentally investigates the convergence of the various schemes when using smooth functions.
      This step tests the ability of the different methods to accurately represent smooth functions as the resolution increases. 
      For the Shape-Preserving Spline (SPS) \cite{10.1145/264029.264050, 10.1145/264029.264059}, DBI and PPI methods, we also investigate the approximation accuracy obtained with varying interpolant polynomial degrees.
    \item The third step focuses on the ability of the different methods to represent a set of challenging test functions with large gradients and/or discontinuities.
This step represents situations often encountered in computational science problems, such as mapping between physics and dynamics meshes in NEPTUNE. 
  \end{itemize}

\section{Positivity-Preserving Interpolants}
\label{sec:positvity}
Preserving positivity while maintaining accuracy is perhaps the key property needed when mapping from one mesh to another in NEPTUNE and 
similar applications. 
This section compares the PCHIP, MQS, SPS, DBI, and PPI against a standard interpolation method using five examples.
The standard polynomial interpolation approach (STD) uses the points in each element to build a standard Lagrange interpolant for that element.
In each of the examples, the different interpolants are constructed using:
  \begin{enumerate}
  \item a uniform mesh that is constructed using uniformly spaced points. 
        In this mesh, all the elements have the same size and the nodes are uniformly spaced inside each element.
  \item an LGL mesh that consists of uniform elements with Legendre Gauss-Lobatto (LGL) quadrature nodes ~\cite{hale2013fast}. 
       For example, the global interval is divided into $(N_{i}-1)/j$ and $(N_{i}-1)^{2} /j^{2}$ elements for the 1D and 2D examples, respectively.
      $N_{i}$ is the total number of points used for the approximation and $j$ is the target polynomial degree.
      For the 1D examples, $j+1$ LGL quadrature nodes are placed inside each element.
      The 2D examples use a tensor product of $(j+1)\times (j+1) =(j+1)^{2}$ LGL nodes in each element.
      This LGL mesh is different than the one used in ~\cite{ouermi2022eno} where the node distribution in each element remains fixed ($(N_{i}-1)/8$ and $(N_{i}-1)^{2}/ 16$) as the target polynomial changes.
      This fixed number of nodes in each element enables the study of the DBI and PPI methods in the case of a fixed LGL mesh with varying polynomial degrees.
  \end{enumerate}
In the figures presented in this section, the black and red plots represent the underlying function and its approximation using the different interpolation methods.   
Both the DBI and PPI methods use a mesh point selection method that favors a symmetric stencil about $x_{i}$ by setting $st=1$ with $\epsilon_{0}=0.01$ and $\epsilon_{1}=1.0$.
The results in Figures \ref{fig:Runge1} to \ref{fig:GelbT2} below demonstrate that the PCHIP, MQS, DBI, SPS, and PPI methods preserve positivity, whereas the standard interpolation methods lead to oscillations and fail to preserve positivity.
Using an LGL mesh reduces oscillations compared to the uniform mesh, but does not guarantee that the interpolating polynomials will be positive.

\subsection{Example I $f_1(x)$}
\label{subsec:example1}
\begin{figure}[H]
    \centering
    \includegraphics[scale=0.5]{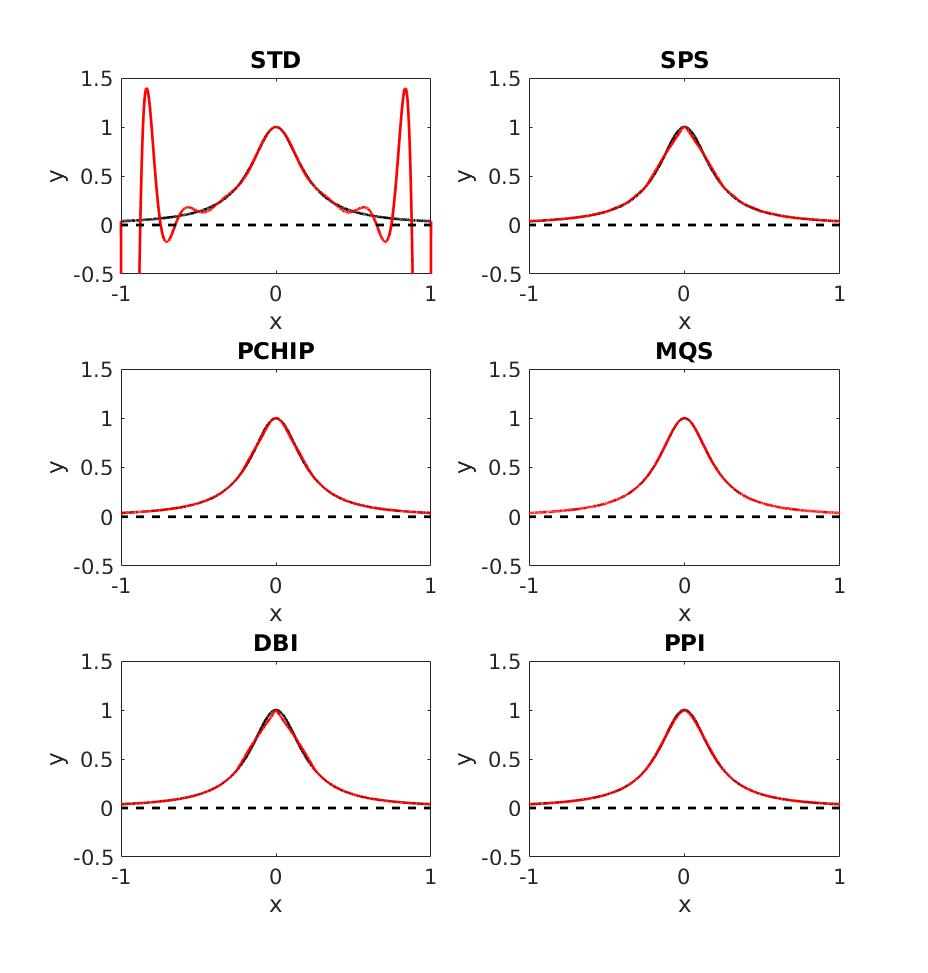}
    \caption{Approximation of the Runge function with the $N=17$ points that are uniformly distributed on the interval $[-1,1]$.}
    \label{fig:Runge1}
  \end{figure}

\begin{figure}[H]
    \centering
    \includegraphics[scale=0.50]{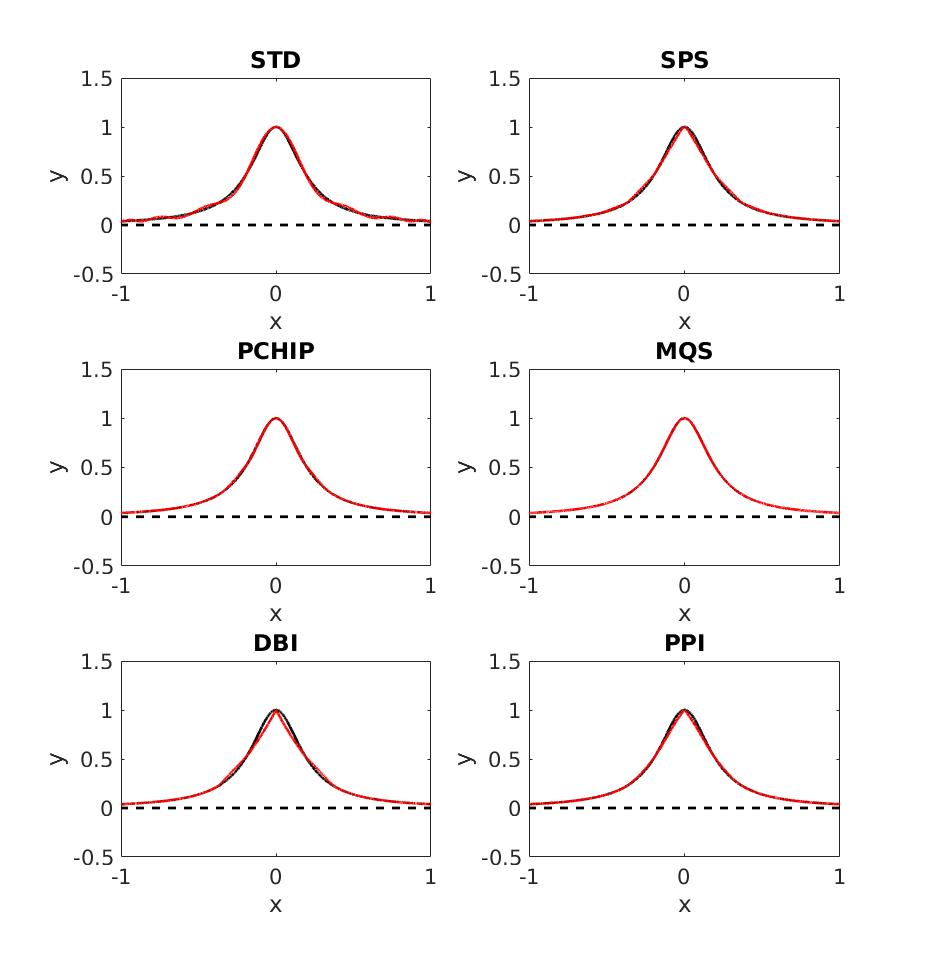}
    \caption{Approximation of the Runge function with the $N=17$ LGL quadrature points distributed on the interval $[-1,1]$.}
    \label{fig:Runge2}
  \end{figure}

This example uses the famous Runge function \cite{doi:10.1080/00029890.1987.12000642} defined as follows:
  \begin{equation}\label{eq:f1}
     f_{1}(x) = \frac{1}{1+25x^{2}}, \quad x \in [-1,1].
  \end{equation}
  Figures \ref{fig:Runge1} and \ref{fig:Runge2} show the different polynomial approximations for this function using $17$ uniformly spaced and LGL points, respectively.
  The target polynomial degree for the standard interpolation, DBI, and PPI is set to $d=16$.
  The standard polynomial interpolation approach, STD, does not preserve positivity with the uniform mesh and generates oscillations in both meshes.
  The PCHIP, MQS, DBI, SPS and PPI methods preserve positivity for both the uniform and LGL meshes.
\subsection{Example II $f_2(x)$}
\label{subsec:examples}
\begin{figure}[H]
    \centering
    \includegraphics[scale=0.50]{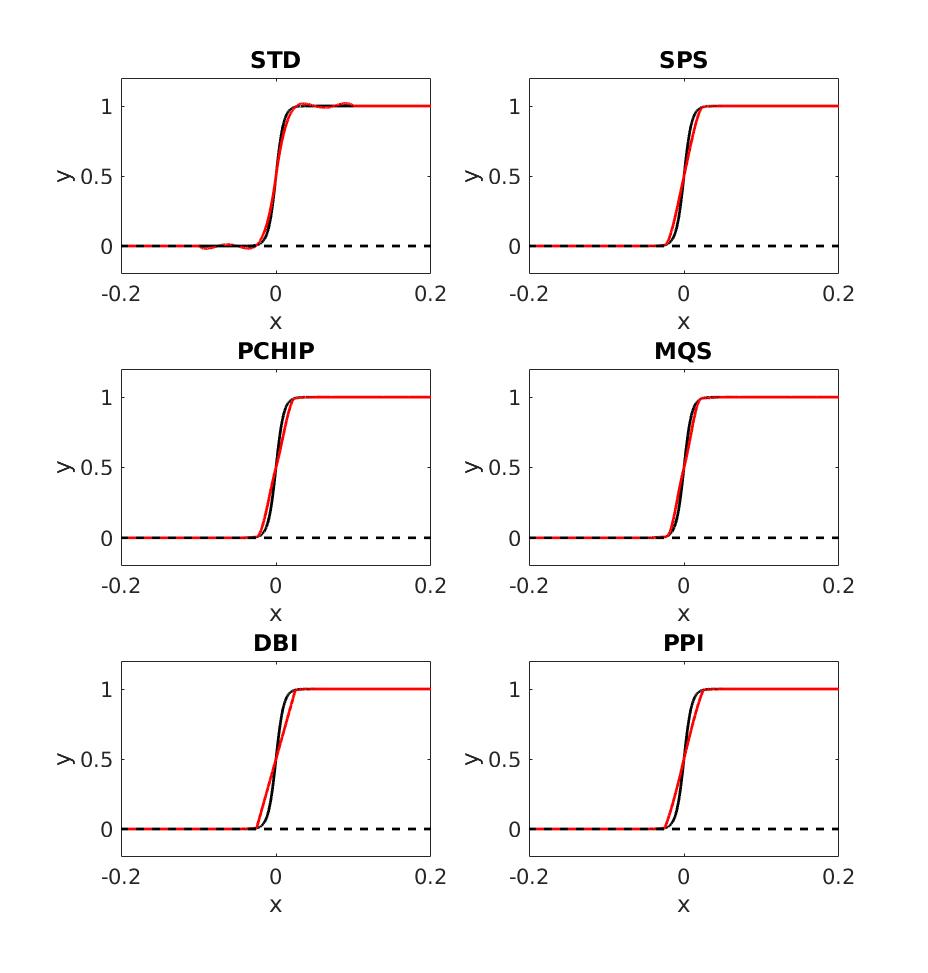}
    \caption{Approximation of $f_{2}(x) = \frac{1}{1+ e^{-2kx}}, \quad k=100 \textrm{, and } x \in [-0.2,0.2]$, with $N=17$ points.
      The points are uniformly distributed and the target polynomial degree for the DBI and PPI is $d=8$.}
    \label{fig:Heaviside1}
  \end{figure}

\begin{figure}[H]
    \centering
    \includegraphics[scale=0.50]{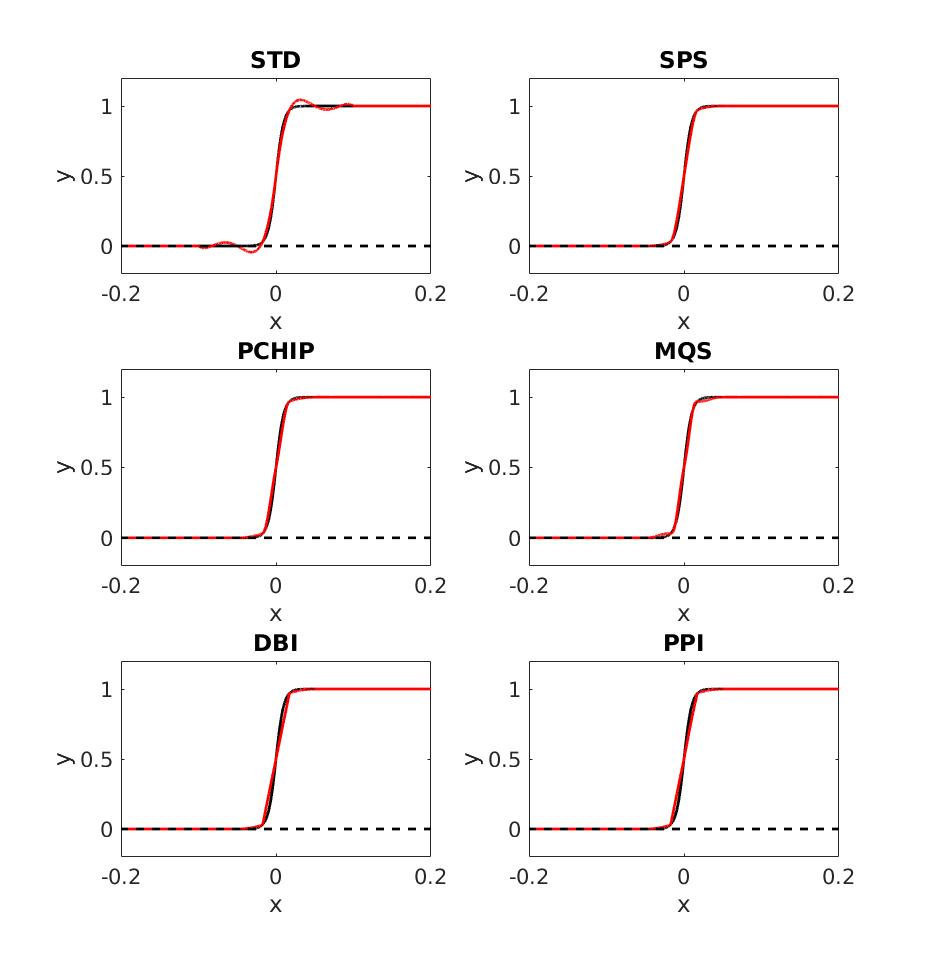}
    \caption{Approximation of $f_{2}(x) = \frac{1}{1+ e^{-2kx}}, \quad k=100 \textrm{, and } x \in [-0.2,0.2]$, with $N=17$ points.
      The interval $[-0.2, 0.2]$ is divided in two elements and $9$ LGL quadrature points are used in each interval.}
    \label{fig:Heaviside2}
  \end{figure}

The second example uses an analytic approximation of the Heaviside function defined as follows:
\begin{equation}\label{eq:f2}
     f_{2}(x) = \frac{1}{1+ e^{-2kx}}, \quad k=100 \textrm{, and } x \in [-0.2,0.2].
\end{equation}
A polynomial approximation of $f_{2}(x)$ is challenging because of the large gradient at about $x =0.$
Attempts to use a global polynomial approximation for this function result in unacceptable oscillations and negative values as observed 
in the Runge example above.
Figures \ref{fig:Heaviside1} and \ref{fig:Heaviside2} show interpolations of $f_{2}(x)$ using a uniform mesh of $17$ points and an LGL mesh with two elements 
each with nine LGL quadrature points in each element.
Standard polynomial interpolation, DBI, and PPI are used with an interpolant of degree $d=8$ for each interval.
Standard polynomial interpolation fails to preserve positivity in both uniform and LGL meshes.
The results demonstrate that the PCHIP, MQS, DBI, SPS  and PPI methods preserve positivity with both the uniform and LGL meshes.

\subsection{Example III $f_3(x)$}
\label{subsec:example3}
\begin{figure}[H]
    \centering
    \includegraphics[scale=0.50]{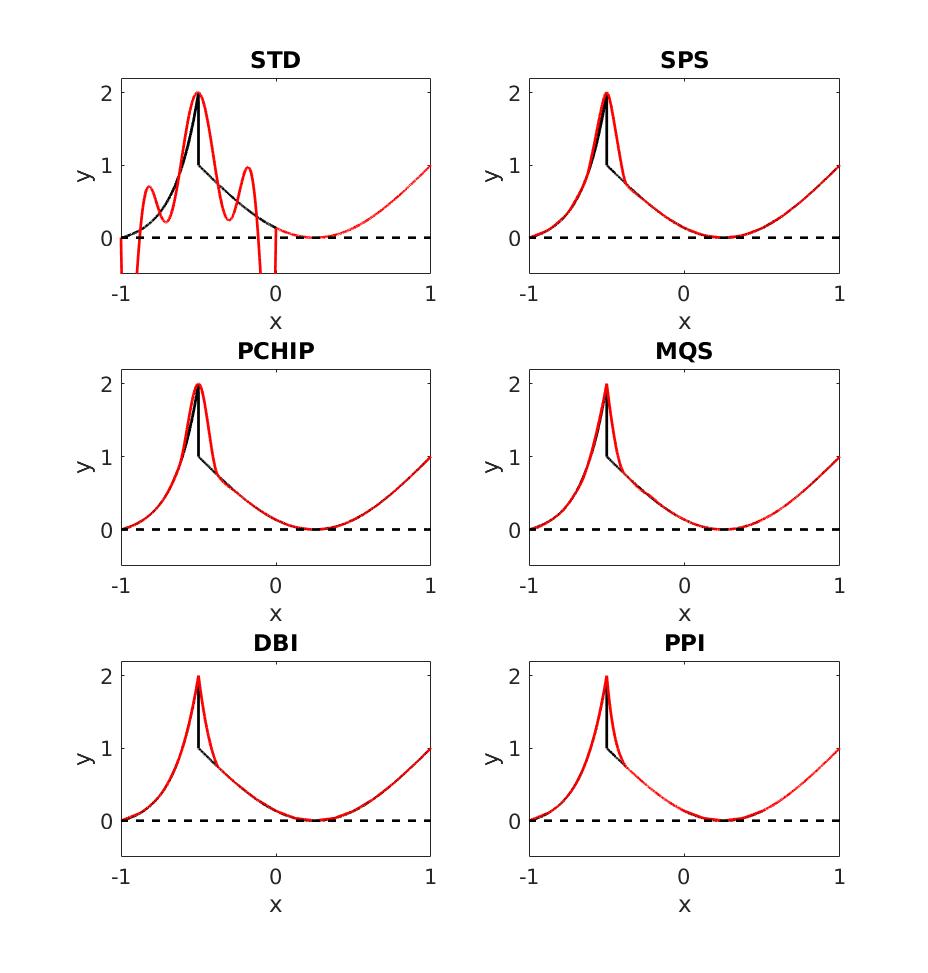}
    \caption{Approximation of $f_{3}(x)$ with $N=17$ points. 
    The points are distributed uniformly over the interval $[-1, 1]$.}
    \label{fig:GelbT1}
  \end{figure}

\begin{figure}[H]
    \centering
    \includegraphics[scale=0.50]{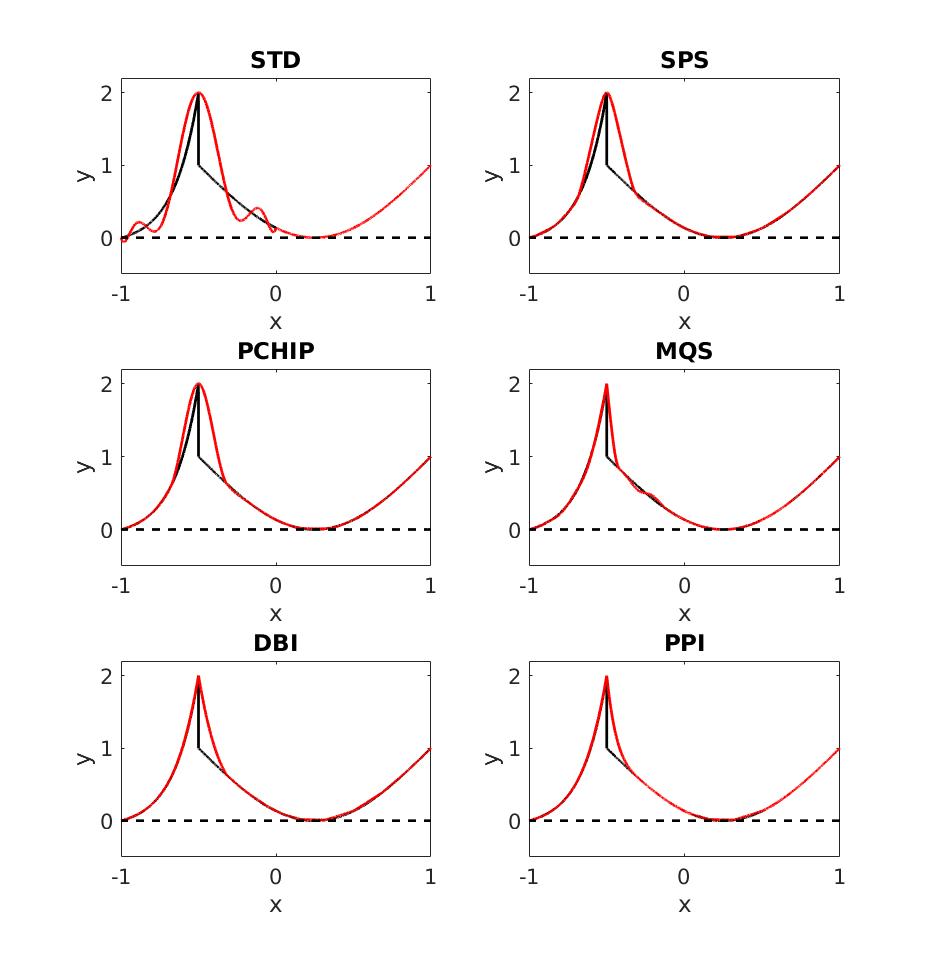}
    \caption{Approximation of $f_{3}(x)$ with $N=17$ points. 
The interval $[-1, 1]$ is divided into four elements, and $d+1$ LGL quadrature points are used in each element.}
    \label{fig:GelbT2}
  \end{figure}

The third example uses a modified version of a function introduced by Tadmor and Tanner \cite{Tadmor2002} and used by Berzins \cite{Berzins} 
in the context of DBI based upon uniform mesh points.
The original function was modified by adding the value one to ensure that the function is positive over the interval $[-1, 1]$.
The modified function is defined as 
  \begin{equation}\label{eq:GT}
    f_{3}(x) =
    \begin{cases}
      1+\frac{2 e^{2\pi (x+1)} -1 -e^{\pi}}{e^{\pi}-1}, \quad x \in [-1, -0.5)\\\\
      1- sin \big(\frac{2\pi x}{3}+ \frac{\pi}{3}\big), \quad x \in [-0.5, 1].
    \end{cases}
  \end{equation}
This function is particularly challenging because of the discontinuity at $x= -0.5$. 
This example uses uniform and LGL meshes of $17$ points.
The LGL mesh consists of four elements and LGL quadrature points are used as the mesh points within each element.
The target interpolant degree for the standard interpolation, DBI, and PPI methods is $d=4$.
Figures \ref{fig:GelbT1} and \ref{fig:GelbT2} demonstrate that the interpolants built using the PCHIP, MQS, SPS, SPS and PPI methods 
remain positive whereas the standard polynomial interpolation approach fails to preserve positivity.
In addition, the oscillations observed with the standard polynomial interpolation method are more pronounced with the uniform mesh 
compared to the LGL mesh. 

\subsection{Example IV $f_4(x)$}
\label{subsec:example4}
\begin{figure}[H]
    \centering
    \includegraphics[scale=0.50]{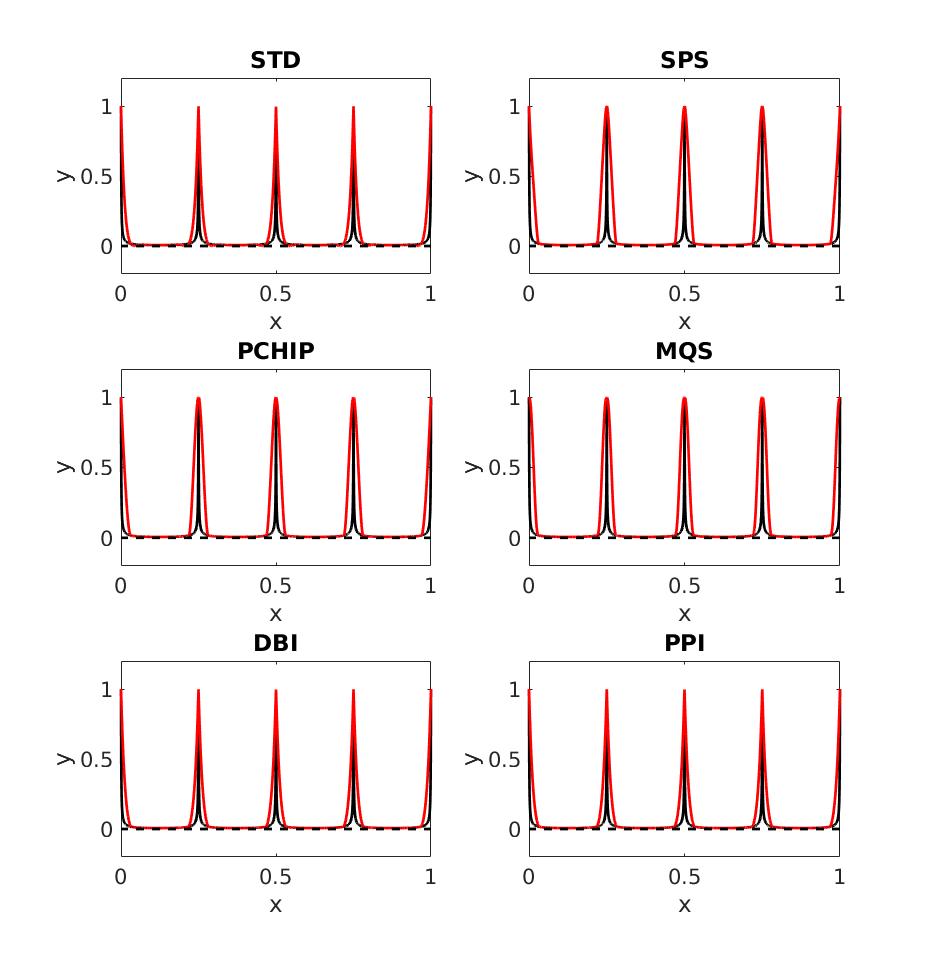}
    \caption{Approximation of $f_{4}(x)$, with $N=33$ points.
      The points are uniformly distributed, and the target polynomial degree for the DBI and PPI is $d=8$.}
    \label{fig:Square1}
  \end{figure}

\begin{figure}[H]
    \centering
    \includegraphics[scale=0.50]{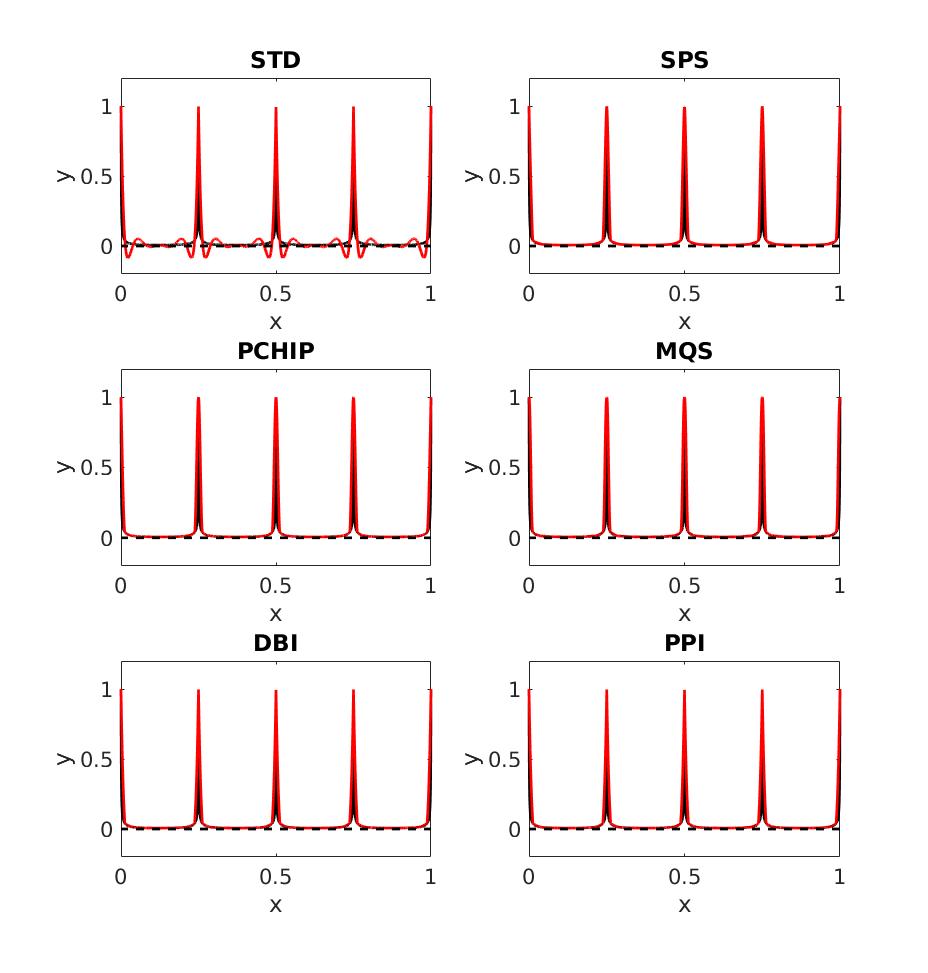}
    \caption{Approximation of $f_{4}(x)$, with $N=33$ points.
      The interval $[0, 1]$ is divided into four elements, and $9$ quadrature points are used in each interval.}
    \label{fig:Square2}
  \end{figure}

  This example consists of a function with multiple spikes defined as follows:
  \begin{equation}\label{eq:square}
    f_{4}(x) = 1.0 - \bigg| \frac{2}{\pi} arctan \bigg( \frac{sin\big( \pi \frac{x}{h}\big)}{\delta}\bigg)\bigg|,\quad x \in [0,1], 
  \end{equation}
  where $h$ represent the element size, and $\delta = 0.01$.
  $f_{4}(x)$ depends on the element size $h$ and therefore, on the number of element in a given interval.
  At the element boundaries, $f_{4}(x)$ is $C^{0}$-continuous with large gradients of opposite signs.  
  This example uses $33$ points, four elements, and nine points in each element.
  The approximations in Figures \ref{fig:Square1} and \ref{fig:Square2} use uniform and LGL quadrature points, respectively. 
  The plots in Figures \ref{fig:Square1} and \ref{fig:Square2} show the standard polynomial interpolation approach lead to oscillation and negative values, whereas the PCHIP, MQS, SPS, DBI, and PPI methods preserve positivity and remove the oscillations.


\subsection{Example V $f_5(x)$}
\label{subsec:example5}
\begin{figure}[H]
    \centering
    \includegraphics[scale=0.50]{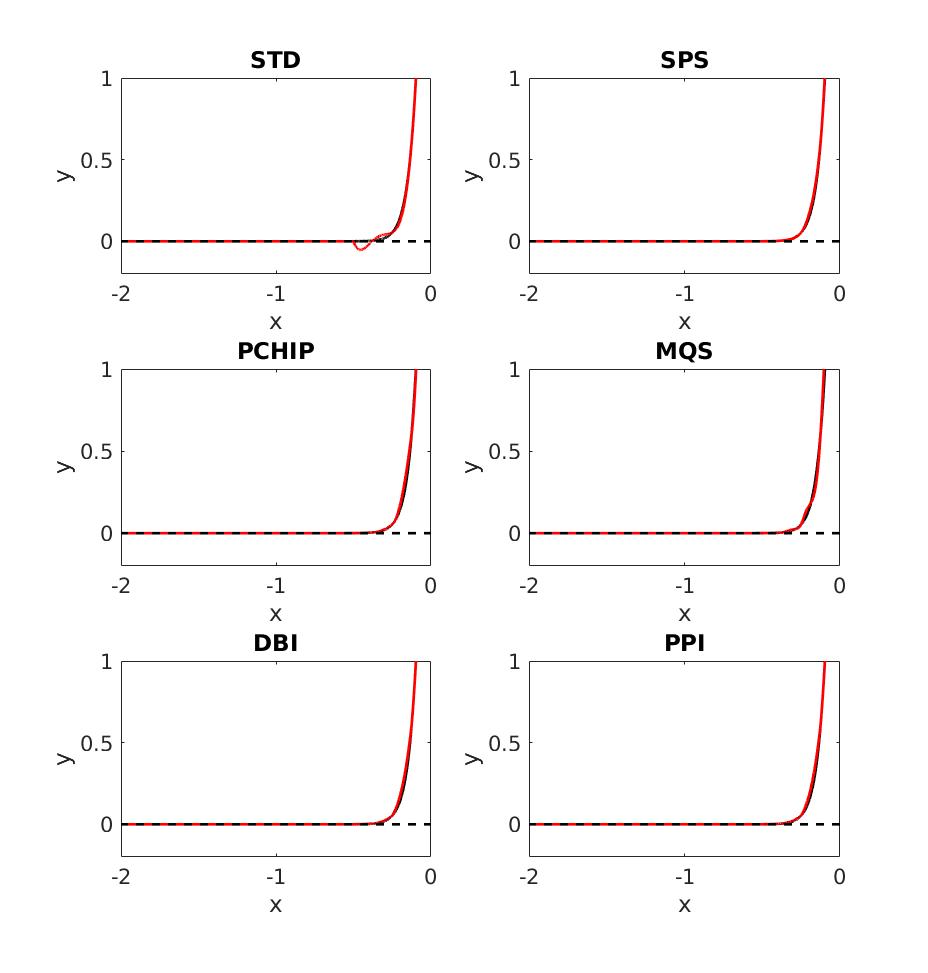}
    \caption{Approximation of $f_{5}(x)$, with $N=17$ points.
      The points are uniformly distributed and the target polynomial degree for the DBI and PPI is $d=4$.}
    \label{fig:Tanh1}
  \end{figure}

\begin{figure}[H]
    \centering
    \includegraphics[scale=0.50]{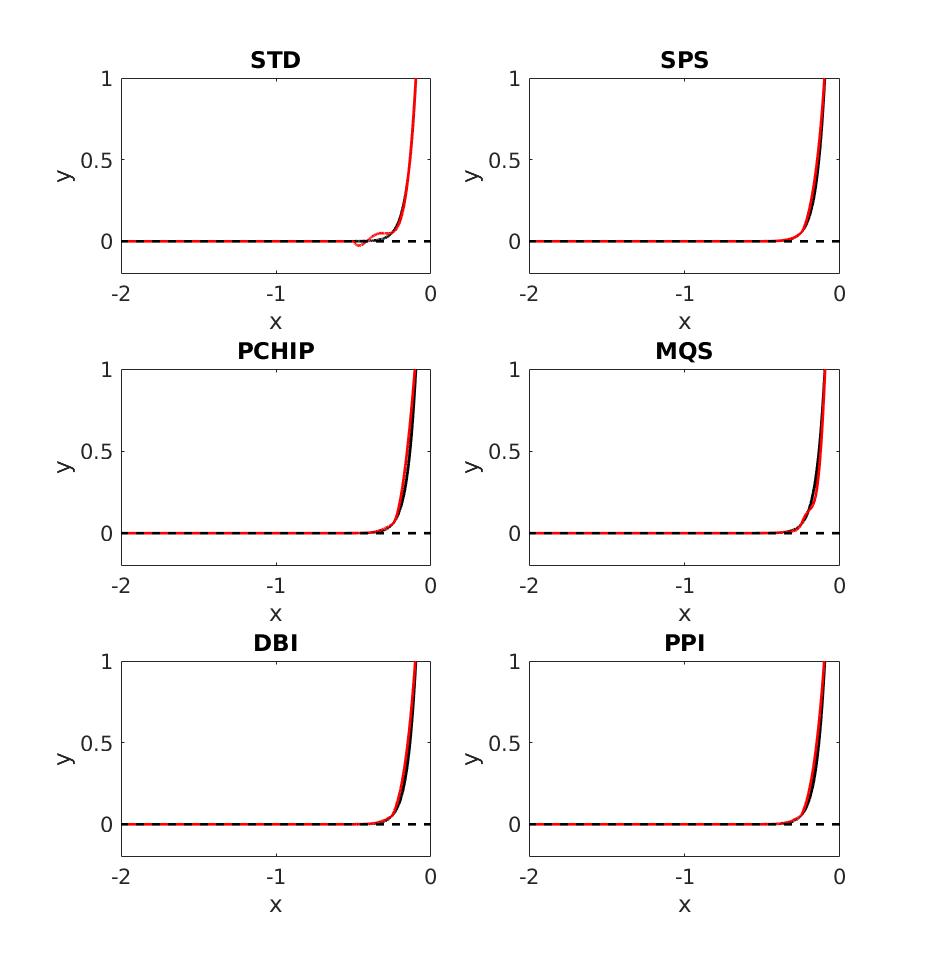}
    \caption{Approximation of $f_{5}(x)$, with $N=17$ points.
      The interval $[-2, 0]$ is divided in four elements and $5$ quadrature points are used in each interval.}
    \label{fig:Tanh2}
  \end{figure}

  This example is constructed using the $tanh$ function and by introducing $C^{0}$-continuities at the elements boundaries.
  The constructed function is defined as follows:
  \begin{equation}\label{eq:tanh}
    f_{5}(x) = 
    \begin{cases}
      tanh(xk)                 			& \textrm{if } x \in [a, a+h]\\
      2tanh(xk) - tanh((a+h)k)  	        & \textrm{if } x \in [a+h, a+2h]\\
      3tanh(xk) - tanh((a+h)k)  - tanh((a+2h)k)	& \textrm{if } x \in [a+2h, a+3h]\\
      \quad \quad \vdots
    \end{cases}
  \end{equation}
  where the overall interval is $[-2,0]$ with $a=-2$ and $k=10$. $h$ represents the size of each element.
  $f_{5}(x)$ depends on the element size $h$ and, therefore, on the number of elements in a given interval.
  This example is built to mirror the $C^{0}$-continuity at the elements boundaries in the spectral element method used in NEPTUNE.
In this example, the gradients at the elements boundaries are always positive, and are not as large as the ones in $f_{4}(x)$ from Example IV.
The approximations shown in Figures \ref{fig:Tanh1} and \ref{fig:Tanh2} use $17$ points, four elements,  and five points inside each element.
  The plots in Figures \ref{fig:Tanh1} and \ref{fig:Tanh2} show that the standard interpolation method does not preserve positivity and that the PCHIP, MQS, SPS, DBI, and PPI can be used to enforce positivity as required.  

\section{Convergence}
  \label{sec:convergence}

This section focuses on the second comparison criterion, which consists of evaluating the convergence of the different methods when applied  to a smooth function.
As NUMA \cite{NUMA}, the dynamics part of NEPTUNE uses a spectral element method that has high-order accuracy, especially in smooth regions.
It is important when interpolating solution values between dynamics and physics meshes for the interpolation scheme to not degrade the accuracy 
obtained from the spectral element method. 
 
The test function 
  \begin{equation}
     f_{6}(x) = 1+sin(x), x \in [0, \pi]
  \end{equation} 
is used to study the convergence of the different methods.
$f_{6}(x)$ is infinitely smooth with no sharp gradients or discontinuities.
These characteristics make $f_{6}(x)$ a suitable test function for evaluating which approach is a good choice for representing smooth functions. 
These experiments focus on the accuracy of the approximation as the resolution and the polynomial degree both increase.

Table \ref{tab:sin_1} shows $L^{2}$-errors when approximating $f_{6}(x)$ using the different interpolation methods.
In this experiment, the parameters $\epsilon_{0}$, $\epsilon_{1}$, and $st$ are chosen to be $0.01$, $1$, and $1$, respectively.
In all cases, the $L^{2}$-error is estimated by sampling the error at $10000$ equally spaced points in the interval and using trapezoidal quadrature.

  Table \ref{tab:rate_sin_1} shows the ratio, $e_{N_{i}} /e_{N_{i+1}}$ of the $L^{2}$-errors in Table \ref{tab:sin_1} as the resolution increases. 
  The DBI and PPI methods lead to smaller error compared to the PCHIP, MQS, and SPS methods.
  As the average polynomial degree increases, the approximation using the DBI method does not improve because the global error is dominated by the local error from the intervals using lower degree interpolants compared to PPI.
  These results show that the conditions for data-boundedness may be more restrictive when it comes to enforcing positivity.
  The SPS method shows smaller errors compared to the other methods.
  Furthermore, as the polynomial degree increases, the accuracy of the approximation decreases.
  These results are consistent with those in \cite{10.1145/264029.264050, 10.1145/264029.264059}.
  Costantini \cite{Costantini1990, 10.1145/264029.264050} demonstrated that the SPS method is bounded by $O(h^{4})$ and in the limit (as the spline degree increases) the spline tends to a linear interpolation.
  The B-spline and PPI methods have smaller $L^{2}$-errors compared to the other methods, and their accuracy improves as the polynomial degree increases.
  Table \ref{tab:rate_sin_1} shows that both methods have better convergence rates compared to PCHIP, MQS, and SPS.
  The PPI method leads to slightly smaller errors compared to the unmodified B-spline approach.
  For $\mathcal{P}_{8}$ and $\mathcal{P}_{16}$ the approximation errors are close to machine precision, which explains the slow rate of convergence observed for B-spline and PPI in Table \ref{tab:rate_sin_1}.

  \begin{table}[H]
     \centering
     \scalebox{1.0}{
     \begin{tabular}{| c| c| c| c| c| c c| c c|}
          \hline
        $N_{i}$   & PCHIP          & MQS           &  SPS            & B-spline      &\multicolumn{2}{c|}{DBI} & \multicolumn{2}{c|}{PPI}     \\
          \hline                         
                  & $L^{2}$-error  & $L^{2}$-error & $L^{2}$-error & $L^{2}$-error & $L^{2}$-error & avg. deg. & $L^{2}$-error & avg. deg. \\
          \hline                           
                  &           &          &         &         & \multicolumn{4}{c|}{$\mathcal{P}_{1}$} \\
          \hline                                               
        $17$      & --        & --       & --      & --      & 2.49E-3 & 1    & 2.49E-3  & 1  \\
        $33$      & --        & --       & --      & --      & 6.22E-4 & 1    & 6.22E-4  & 1  \\
        $65$      & --        & --       & --      & --      & 1.56E-4 & 1    & 1.56E-4  & 1  \\
        $129$     & --        & --       & --      & --      & 3.89E-5 & 1    & 3.89E-5  & 1  \\
        $257$     & --        & --       & --      & --      & 9.72E-6 & 1    & 9.72E-6  & 1  \\
          \hline                         
                  & $\mathcal{P}_{3}$ & $\mathcal{P}_{5}$      & \multicolumn{6}{c|}{$\mathcal{P}_{4}$} \\
          \hline                         
        $17$      & 4.49E-4   & 4.47E-5  & 4.84E-4 & 4.52E-6   &  6.70E-06   & 3.94   & 2.52E-06   & 4  \\
        $33$      & 7.83E-5   & 7.42E-6  & 1.20E-4 & 2.07E-7   &  2.85E-07   & 3.97   & 6.94E-08   & 4  \\
        $65$      & 1.38E-5   & 1.31E-6  & 3.01E-5 & 1.22E-8   &  1.24E-08   & 3.98   & 1.96E-09   & 4  \\
        $129$     & 2.45E-6   & 2.31E-7  & 7.52E-6 & 7.56E-10  &  5.43E-10   & 3.99   & 5.73E-11   & 4  \\
        $257$     & 4.34E-7   & 4.09E-8  & 1.88E-6 & 4.72E-11  &  2.39E-11   & 4.00   & 1.73E-12   & 4  \\
          \hline                                                                                       
                  &           &          & \multicolumn{6}{c|}{$\mathcal{P}_{8}$} \\
          \hline                                                 
                                                                                       
        $17$      & --        & --       & 2.00E-3 & 2.45E-9   &  6.21E-06   & 7.69   & 1.06E-09   & 8  \\
        $33$      & --        & --       & 4.96E-4 & 3.47E-12  &  2.76E-07   & 7.84   & 1.83E-12   & 8  \\
        $65$      & --        & --       & 1.24E-4 & 6.11E-15  &  1.22E-08   & 7.92   & 3.44E-15   & 8  \\
        $129$     & --        & --       & 3.10E-5 & 3.23E-15  &  5.40E-10   & 7.96   & 1.00E-15   & 8  \\
        $257$     & --        & --       & 7.74E-6 & 2.93E-15  &  2.39E-11   & 7.98   & 9.64E-16   & 8  \\
          \hline                                                                                          
                  &           &          & \multicolumn{6}{c|}{$\mathcal{P}_{16}$} \\
          \hline                                                 
        $17$      & --        & --       & 3.10E-3 & 5.61E-15  &  6.21E-06   & 15.19   & 3.98E-15  & 16  \\
        $33$      & --        & --       & 7.75E-4 & 4.35E-13  &  2.76E-07   & 15.59   & 1.95E-15  & 16  \\
        $65$      & --        & --       & 1.94E-4 & 2.75E-13  &  1.22E-08   & 15.80   & 4.65E-15  & 16  \\
        $129$     & --        & --       & 4.84E-5 & 9.00E-14  &  5.40E-10   & 15.90   & 2.33E-15  & 16  \\
        $257$     & --        & --       & 1.21E-5 & 6.83E-14  &  2.39E-11   & 15.95   & 1.10E-15  & 16  \\
         \hline              
     \end{tabular}}                                                                     
     \caption{$L^2$-errors when using the PCHIP, MQS,  SPS, B-splines, DBI, and PPI methods to approximate the function $f_{6}(x)$.
           $N_{i}$ represents the number of input points used to build the approximation.
           $P_{j}$ represents the space of polynomials of degree $j$, with $j$ being the target degree for each interval.
           The seventh and ninth columns show the average polynomial degree used for the DBI and PPI methods, respectively.
           The input points are uniformly distributed over the interval $[0, \pi]$.} 
     \label{tab:sin_1}
   \end{table}

  \begin{table}[H]
     \centering
     \scalebox{1.0}{
     \begin{tabular}{| c| c  c  c  c  c  c|}
          \hline
        $e_{N_{i}}/e_{N_{i+1}}$    & PCHIP  & MQS   & SPS  & B-spline   & DBI  & PPI     \\
         \hline                         
                          &        &       &    &   & \multicolumn{2}{|c|}{$\mathcal{P}_{1}$} \\
          \hline                                               
        $e_{17}/e_{33}$   & --     &  --   & --     & --      &   4  &  4    \\
        $e_{33}/e_{65}$   & --     &  --   & --     & --      &   4  &  4    \\
        $e_{65}/e_{129}$  & --     &  --   & --     & --      &   4  &  4    \\
        $e_{129}/e_{257}$ & --     &  --   & --     & --      &   4  &  4    \\
          \hline                                                          
                          & $\mathcal{P}_{3}$ & $\mathcal{P}_{5}$ &  \multicolumn{4}{|c|}{$\mathcal{P}_{4}$} \\
          \hline                                                    
        $e_{17}/e_{33}$   & 5.73   & 6.02  & 4.03 & 21.84      &  24  & 36    \\
        $e_{33}/e_{65}$   & 5.67   & 5.67  & 3.99 & 16.97      &  23  & 35    \\
        $e_{65}/e_{129}$  & 5.63   & 5.66  & 4.00 & 16.13      &  23  & 34    \\
        $e_{129}/e_{257}$ & 5.64   & 5.66  & 4.00 & 16.02      &  23  & 33    \\
          \hline                                                           
                  &       &        &  \multicolumn{4}{|c|}{$\mathcal{P}_{8}$} \\
          \hline                                                 
        $e_{17}/e_{33}$   & --     & --    & 4.03 & 706.05     &  22  & 576  \\
        $e_{33}/e_{65}$   & --     & --    & 4.00 & 567.92     &  23  & 533  \\
        $e_{65}/e_{129}$  & --     & --    & 4.00 & 1.89       &  23  &  3   \\
        $e_{129}/e_{257}$ & --     & --    & 4.01 & 1.10       &  23  &  1   \\
          \hline                                                           
                  &       &        &  \multicolumn{4}{|c|}{$\mathcal{P}_{16}$} \\
          \hline                                                 
        $e_{17}/e_{33}$   & --     & --    & 4.01 & 0.01       &  22 &  2    \\
        $e_{33}/e_{65}$   & --     & --    & 3.99 & 1.58       &  23 &  0    \\
        $e_{65}/e_{129}$  & --     & --    & 4.01 & 3.06       &  23 &  2    \\
        $e_{129}/e_{257}$ & --     & --    & 4.00 & 1.32       &  23 &  2    \\
         \hline              
     \end{tabular}}                                                                     
     \caption{Ratio of $L^2$-errors from Table \ref{tab:sin_1} ($e_{N_{i}}/e_{N_{i+1}}$). 
           $N_{i}$ represents the number of input points used to build the approximation.
           $P_{j}$ represents the space of polynomials of degree $j$, with $j$ being the target degree for each interval.}
     \label{tab:rate_sin_1}
   \end{table}

\section{Results}
\label{sec:results}
  In this section, the different interpolation methods are used to approximate functions with steep gradients, $C^{0}$-continuity, and discontinuities.
  These experiments focus on the third criterion, which consists of  evaluating the ability of the different methods to represent nonsmooth functions.
  The data points for the interpolation are sampled from 1D and 2D functions.
  Two types of meshes are used for the various experiments.
  The first type of mesh uses uniform elements and uniformly spaced nodes within each element.
  The second type of mesh uses uniform elements and Legendre Gauss-Lobatto (LGL) quadrature nodes \cite{hale2013fast} within each element.
  The experimental results compare the DBI and PPI methods against the SPS, PCHIP \cite{fritsch1980monotone}, and MQS \cite{Lux2019ANAF} methods.    

  The MQS algorithm is designed for monotonically increasing data.
  In order to use the MQS approach with the different 1D examples, we divide the data into monotonically increasing and decreasing regions.
  For the monotonically increasing data, the MQS algorithm is applied directly. 
  For the monotonically decreasing data, we uses the reflection of the data about a vertical axis and applied the MQS algorithm.
  Because of the data transformation involved, the MQS method is used only for the 1D examples.
   
 Tables \ref{tab:runge1D_1} -- \ref{tab:GelbT_2} and \ref{tab:runge2D_1} -- \ref{tab:T2_2} show $L^{2}$-errors when using the different methods to approximate the 1D and 2D functions, respectively.
  The `avg. deg.' columns in these  tables represent the average polynomial degree used in the DBI and PPI method.

  \subsection{Example I $f_1(x)$}
  \label{subsec:example1r}
  This example is the 1D Runge function \cite{doi:10.1080/00029890.1987.12000642} defined in Equation \ref{eq:f1} with $\epsilon_{0}=0.01$, $\epsilon_{1}=1$ and $st=2$.
  Tables \ref{tab:runge1D_1} and \ref{tab:runge1D_2} demonstrate that the PPI method gives smaller approximation 
errors when compared to the other approaches.
  The requirements of data boundedness in the DBI method are restrictive compared to the positivity requirements in PPI.
  These restrictions lead to lower average polynomial degrees for DBI compared to PPI, as shown in the sixth and eighth columns in Tables \ref{tab:runge1D_1} and \ref{tab:runge1D_2}.
  In the case of the uniform mesh, as the average degree used by DBI increases the $L^{2}$-errors remain the same.
  The approximation error using the DBI method does not improve as the average degree increases because the global error is dominated by the local error of those subintervals with low degree interpolants.
  The polynomial degree of the interpolants used for these intervals remains the same as the average polynomial degree of the interpolant increases elsewhere.
  The PPI methods uses higher order interpolants compared to the SPS, DBI, PCHIP, and MQS methods in both the uniform and LGL meshes. 
  The uniform mesh leads to slightly more accurate results than the LGL mesh.
  These results show that the PPI method is a suitable approach for interpolating data from one mesh to another in cases where the underlying 
 function is similar to the Runge function.  

    \begin{table}[H]
     \centering
     \scalebox{1.0}{
     \begin{tabular}{| c| c| c| c| c c| c c|}
          \hline
        $N_{i}$   & PCHIP     & MQS      & SPS & \multicolumn{2}{c|}{DBI} & \multicolumn{2}{c|}{PPI}     \\
          \hline                           
                  & $L^{2}$-error & $L^{2}$-error & $L^{2}$-error & $L^{2}$-error & avg. deg. & $L^{2}$-error & avg. deg. \\
          \hline                           
                  &           &          &          & \multicolumn{4}{c|}{$\mathcal{P}_{1}$} \\
          \hline                           
        $17$      & --        & --       & --       & 2.16E-2   & 1      & 2.16E-2    & 1   \\
        $33$      & --        & --       & --       & 6.02E-3   & 1      & 6.02E-3    & 1   \\
        $65$      & --        & --       & --       & 1.52E-3   & 1      & 1.52E-3    & 1   \\
        $129$     & --        & --       & --       & 3.82E-4   & 1      & 3.82E-4    & 1   \\
        $257$     & --        & --       & --       & 9.56E-5   & 1      & 9.56E-5    & 1   \\
          \hline                                                                               
                  & $\mathcal{P}_{3}$  &  $\mathcal{P}_{5}$   & \multicolumn{5}{c|}{$\mathcal{P}_{4}$} \\
          \hline                                      
        $17$      & 7.15E-3   & 5.72E-3  &  8.34E-3 & 8.34E-3   & 4      & 7.02E-3   & 4    \\
        $33$      & 1.91E-3   & 3.95E-4  &  5.91E-4 & 5.91E-4   & 4      & 5.91E-4   & 4    \\
        $65$      & 3.70E-4   & 6.44E-5  &  4.26E-5 & 4.26E-5   & 3.98   & 2.39E-5   & 4    \\
        $129$     & 6.79E-5   & 5.27E-6  &  2.68E-6 & 2.68E-6   & 3.98   & 8.00E-7   & 4    \\
        $257$     & 1.22E-5   & 6.83E-7  &  8.63E-8 & 8.63E-8   & 4.00   & 2.55E-8   & 4    \\
          \hline                                                                                  
                  &           &          &         \multicolumn{5}{c|}{$\mathcal{P}_{8}$} \\
          \hline                                         
        $17$      & --        & --       & 1.21E-2  & 4.61E-3   & 7.88   & 3.11E-3   & 7.88  \\ 
        $33$      & --        & --       & 2.74E-3  & 4.43E-4   & 7.88   & 1.51E-4   & 8     \\ 
        $65$      & --        & --       & 6.86E-4  & 3.67E-5   & 7.92   & 1.05E-6   & 8     \\ 
        $129$     & --        & --       & 1.72E-4  & 2.56E-6   & 7.92   & 3.10E-9   & 8     \\ 
        $257$     & --        & --       & 4.30E-5  & 8.24E-8   & 7.97   & 6.80E-12  & 8     \\ 	  
          \hline                                         
                  &           &          &         \multicolumn{5}{c|}{$\mathcal{P}_{16}$} \\
          \hline                                         
	$17$      & --        & --       & 1.64E-2  & 4.34E-3   & 11.31   & 3.44E-3   & 11.75    \\  
	$33$      & --        & --       & 4.25E-3  & 4.21E-4   & 15.62   & 4.85E-5   & 16    \\ 
	$65$      & --        & --       & 1.07E-3  & 3.67E-5   & 15.69   & 5.92E-8   & 16    \\
	$129$     & --        & --       & 2.69E-4  & 2.56E-6   & 15.80   & 4.21E-12  & 16    \\
	$257$     & --        & --       & 6.71E-5  & 8.24E-8   & 15.91   & 2.18E-16  & 16    \\
          \hline              
     \end{tabular}}                                                                     
     \caption{$L^2$-errors when using the PCHIP, MQS, SPS, DBI, and PPI methods to approximate the Runge function $f_{1}(x) = \frac{1}{1+25x^{2}}, \quad x \in [-1,1]$.
           $N_{i}$ represents the number of input points used to build the approximation.
           $P_{j}$ represents the use of polynomials of degree $j$, with $j$ being the target degree for each interval.
           The sixth and eighth columns show the average polynomial degree used for the DBI and PPI methods, respectively.
           The input points are uniformly distributed over the interval $[-1,1]$.} 
     \label{tab:runge1D_1}
   \end{table}
   \begin{table}[H]
     \centering
     \scalebox{1.0}{
     \begin{tabular}{| c| c| c| c| c c| c c|}
          \hline
        $N_{i}$   & PCHIP     & MQS      &SPS      & \multicolumn{2}{c|}{DBI} & \multicolumn{2}{c|}{PPI}     \\
          \hline                          
                  & $L^{2}$-error & $L^{2}$-error & $L^{2}$-error & $L^{2}$-error & avg. deg. & $L^{2}$-error & avg. deg. \\
          \hline                           
                  &           &          &          & \multicolumn{4}{c|}{$\mathcal{P}_{1}$} \\
          \hline                           
        $17$      & --        & --       & --       & 2.16E-2    & 1      & 2.16E-2   & 1  \\
        $33$      & --        & --       & --       & 6.02E-3    & 1      & 6.02E-3   & 1  \\
        $65$      & --        & --       & --       & 1.52E-3    & 1      & 1.52E-3   & 1  \\
        $129$     & --        & --       & --       & 3.82E-4    & 1      & 3.82E-4   & 1  \\
        $257$     & --        & --       & --       & 9.56E-5    & 1      & 9.56E-5   & 1  \\
          \hline                         
                  & $\mathcal{P}_{3}$ & $\mathcal{P}_{5}$ &  \multicolumn{5}{c|}{$\mathcal{P}_{4}$} \\
          \hline                         
        $17$      & 1.02E-2   & 6.29E-3  & 9.73E-3  &  8.63E-3   & 4      & 8.39E-3   & 4    \\
        $33$      & 1.86E-3   & 9.13E-4  & 1.63E-3  &  7.95E-4   & 4      & 7.80E-4   & 4    \\
        $65$      & 3.68E-4   & 8.47E-5  & 2.24E-4  &  4.76E-5   & 3.98   & 4.64E-5   & 4    \\
        $129$     & 7.20E-5   & 6.23E-6  & 6.03E-5  &  1.49E-6   & 3.98   & 1.27E-6   & 4    \\
        $257$     & 1.52E-5   & 5.72E-7  & 1.48E-5  &  4.68E-8   & 4      & 3.95E-8   & 4    \\
          \hline                                                                                 
                  &           &          &        \multicolumn{5}{c|}{$\mathcal{P}_{8}$} \\
          \hline                                      
        $17$      & --        & --       & 8.44E-3  &  3.49E-3   & 8.00   & 4.40E-3    & 8      \\
        $33$      & --        & --       & 2.69E-3  &  1.76E-4   & 7.88   & 1.76E-4    & 8      \\
        $65$      & --        & --       & 7.59E-4  &  3.25E-6   & 7.92   & 3.01E-6    & 8      \\
        $129$     & --        & --       & 2.61E-4  &  5.64E-8   & 7.94   & 8.82E-9    & 8      \\
        $257$     & --        & --       & 6.85E-5  &  3.51E-9   & 7.96   & 3.96E-11   & 8      \\
          \hline                                                                                 
                  &           &          &        \multicolumn{5}{c|}{$\mathcal{P}_{16}$} \\
          \hline                                      
        $17$      & --        & --       & 2.29E-2  &  9.12E-3   & 12.19   & 1.25E-2    & 12.62    \\
        $33$      & --        & --       & 3.26E-3  &  5.87E-5   & 15.28   & 5.86E-5    & 16       \\
        $65$      & --        & --       & 1.11E-3  &  1.41E-7   & 15.62   & 1.17E-7    & 16       \\
        $129$     & --        & --       & 3.12E-4  &  3.52E-9   & 15.90   & 4.44E-11   & 16       \\
        $257$     & --        & --       & 1.08E-4  &  1.56E-10  & 15.91   & 2.88E-15   & 16       \\
          \hline
     \end{tabular}}                                                                                  
     \caption{$L^2$-errors when using the PCHIP, MQS, SPS, DBI, and PPI methods to approximate the Runge function $f_{1}(x) = \frac{1}{1+25x^{2}}, \quad x \in [-1,1]$.
           $N_{i}$ represents the number of input points used to build the approximation.
           $P_{j}$ represents the use of polynomials of degree $j$, with $j$ being the target degree for each interval.
           The sixth and eight columns show the average polynomial degree used for the DBI and PPI methods respectively.
           The interval [-1,1] is divided into $(N_{i}-1) / j$ and $j+1$ LGL quadrature points are used in each element.}
     \label{tab:runge1D_2}
   \end{table}

  \subsection{Example II $f_2(x)$ }
  \label{subsec:example2r}
  The second example uses the analytic approximation of the Heaviside function defined in Equation \ref{eq:f2} with $\epsilon_{0}=0.01$, $\epsilon_{1}=1$ and $st=2$.
  As mentioned in Section 3.2, this function, $f_{2}(x)$, is challenging because of the sharp gradient around $x=0$.
  For polynomial degree five or less, the results from Tables \ref{tab:heaviside_1} and \ref{tab:heaviside_2} suggest 
  that the MQS method leads to slighter better approximations than DBI and PPI for $f_{2}(x)$.
  Overall, the results from Tables \ref{tab:heaviside_1} and \ref{tab:heaviside_2} indicate that the DBI and PPI methods have smaller $L^{2}$-errors compared to the other methods.
  Approximating $f_{2}(x)$ from data on a uniform mesh leads to slightly better results compared to LGL mesh data.
  For smooth data with a large gradient, these results indicate that both the DBI and PPI approaches are suitable for interpolating from one mesh to another.

    \begin{table}[H]
     \centering
     \scalebox{1.0}{
     \begin{tabular}{| c| c| c| c| c c| c c|}
          \hline
        $N_{i}$   & PCHIP     & MQS      & SPS      & \multicolumn{2}{c|}{DBI} & \multicolumn{2}{c|}{PPI}     \\
          \hline                      
                  & $L^{2}$-error & $L^{2}$-error & $L^{2}$-error & $L^{2}$-error & avg. deg. & $L^{2}$-error & avg. deg. \\
          \hline                          
                  &           &          &          & \multicolumn{4}{c|}{$\mathcal{P}_{1}$} \\
          \hline                          
        $17$      & --        & --       & --       & 2.89E-2    &  1     & 2.89E-2  & 1  \\ 
        $33$      & --        & --       & --       & 7.69E-3    &  1     & 7.69E-3  & 1  \\
        $65$      & --        & --       & --       & 1.80E-3    &  1     & 1.80E-3  & 1  \\
        $129$     & --        & --       & --       & 4.58E-4    &  1     & 4.58E-4  & 1  \\
        $257$     & --        & --       & --       & 1.15E-4    &  1     & 1.15E-4  & 1  \\
          \hline                          
                  & $\mathcal{P}_{3}$ & $\mathcal{P}_{5}$ &    \multicolumn{5}{c|}{$\mathcal{P}_{4}$} \\
          \hline                          
        $17$      & 2.02E-02  & 1.67E-2  & 1.82E-2  &  2.23E-2   & 2.75   & 2.23E-2  & 3.38  \\ 
        $33$      & 3.38E-03  &	4.16E-3  & 3.72E-3  &  4.09E-3   & 3.62   & 4.10E-3  & 3.72  \\
        $65$      & 3.59E-04  &	2.29E-4  & 3.40E-4  &  3.05E-4   & 3.86   & 3.05E-4  & 3.86  \\
        $129$     & 4.21E-05  &	7.48E-6  & 5.36E-5  &  1.35E-5   & 3.88   & 1.35E-5  & 3.88  \\
        $257$     & 5.12E-06  &	2.16E-7  & 1.27E-5  &  4.71E-7   & 3.85   & 4.71E-7  & 3.86  \\
          \hline                                                                               
                  &           &          &          \multicolumn{5}{c|}{$\mathcal{P}_{8}$} \\
          \hline                                      
 	$17$      & --        & --       & 3.75E-3  &  2.08E-2   & 3.25   & 2.08E-2   & 5.50  \\
 	$33$      & --        & --       & 5.24E-3  &  3.36E-3   & 3.88   & 3.33E-3   & 5.72  \\
 	$65$      & --        & --       & 8.71E-4  &  1.38E-4   & 7.59   & 1.38E-4   & 7.59  \\
 	$129$     & --        & --       & 2.08E-4  &  1.22E-6   & 7.68   & 1.22E-6   & 7.73  \\
 	$257$     & --        & --       & 5.17E-5  &  4.44E-9   & 7.61   & 4.44E-9   & 7.67  \\
          \hline                                                                               
                  &           &          &          \multicolumn{5}{c|}{$\mathcal{P}_{16}$} \\
          \hline                                      
 	$17$      & --        & --       & 5.90E-3  &  2.00E-2   & 4.25   & 2.00E-2    & 6.62   \\
 	$33$      & --        & --       & 6.34E-3  &  2.93E-3   & 4.38   & 2.91E-3    & 9.72   \\
 	$65$      & --        & --       & 1.30E-3  &  9.17E-5   & 14.64  & 9.17E-5    & 14.86  \\
 	$129$     & --        & --       & 3.23E-4  &  1.70E-7   & 15.15  & 1.70E-7    & 15.41  \\
 	$257$     & --        & --       & 8.08E-5  &  2.64E-11  & 15.05  & 2.64E-11   & 15.30  \\
          \hline
     \end{tabular}}                                                                                  
     \caption{$L^2$-errors when using the PCHIP, MQS, SPS,  DBI, and PPI methods to approximate the function $f_{2}(x) = \frac{1}{1+ e^{-2kx}}, \quad k=100 \textrm{, and } x \in [-0.2,0.2]$.
           $N_{i}$ represents the number of input points used to build the approximation.
           $P_{j}$ represents the use of polynomials of degree $j$, with $j$ being the target degree for each interval.
           The sixth and eighth columns show the average polynomial degree used for the DBI and PPI methods respectively.
           The input points are uniformly distributed over the interval $[-1,1]$.} 
     \label{tab:heaviside_1}
   \end{table}
   \begin{table}[H]
     \centering
     \scalebox{1.0}{
     \begin{tabular}{| c| c| c| c| c c| c c|}
          \hline
        $N_{i}$   & PCHIP     & MQS      & SPS      & \multicolumn{2}{c|}{DBI} & \multicolumn{2}{c|}{PPI}     \\
          \hline                  
                  & $L^{2}$-error & $L^{2}$-error & $L^{2}$-error & $L^{2}$-error & avg. deg. & $L^{2}$-error & avg. deg. \\
          \hline                          
                  &           &          &          & \multicolumn{4}{c|}{$\mathcal{P}_{1}$} \\
          \hline                          
        $17$      &  --        & --       & --       & 2.89E-2  &  1      & 2.89E-2  & 1  \\ 
        $33$      &  --        & --       & --       & 7.69E-3  &  1      & 7.69E-3  & 1  \\
        $65$      &  --        & --       & --       & 1.80E-3  &  1      & 1.80E-3  & 1  \\
        $129$     &  --        & --       & --       & 4.58E-4  &  1      & 4.58E-4  & 1  \\
        $257$     &  --        & --       & --       & 1.15E-4  &  1      & 1.15E-4  & 1  \\
          \hline                  
                  & $\mathcal{P}_{3}$ & $\mathcal{P}_{5}$ & \multicolumn{5}{c|}{$\mathcal{P}_{4}$} \\
          \hline                          
        $17$      & 8.60E-3   & 7.38E-3  & 7.15E-3  &  1.26E-2   & 2.88   & 1.25E-2   & 3.44   \\
        $33$      & 2.50E-3   & 2.50E-3  & 8.04E-4  &  3.11E-3   & 3.03   & 2.83E-3   & 3.44   \\
        $65$      & 6.36E-4   & 2.11E-4  & 4.18E-4  &  3.28E-4   & 3.81   & 3.72E-4   & 3.84   \\
        $129$     & 1.02E-4   & 1.01E-5  & 9.07E-5  &  1.55E-5   & 3.88   & 1.55E-5   & 3.88   \\
        $257$     & 1.83E-5   & 2.93E-7  & 1.83E-5  &  6.29E-7   & 3.85   & 6.29E-7   & 3.86   \\
          \hline                                                                                
                  &           &          &        \multicolumn{5}{c|}{$\mathcal{P}_{8}$} \\
          \hline                                      
        $17$      & --        & --       & 4.43E-3  &  4.87E-3   & 3.50   & 4.68E-3   & 5.00  \\
        $33$      & --        & --       & 2.51E-3  &  8.71E-4   & 4.34   & 7.84E-4   & 5.75  \\
        $65$      & --        & --       & 1.00E-3  &  7.57E-5   & 6.64   & 1.24E-4   & 7.28  \\
        $129$     & --        & --       & 3.65E-4  &  2.17E-6   & 7.65   & 2.17E-6   & 7.73  \\
        $257$     & --        & --       & 9.11E-5  &  1.95E-8   & 7.55   & 1.95E-8   & 7.73  \\
          \hline                                                                                
                  &           &          &        \multicolumn{5}{c|}{$\mathcal{P}_{16}$} \\
          \hline                                      
        $17$      & --        & --       & 4.52E-2  &  3.77E-2   & 3.81    & 3.73E-2    & 7.25   \\
        $33$      & --        & --       & 2.03E-3  &  2.53E-4   & 5.56    & 5.23E-4    & 9.84   \\
        $65$      & --        & --       & 9.55E-4  &  1.37E-5   & 10.53   & 6.95E-5    & 12.56 \\
        $129$     & --        & --       & 4.16E-4  &  2.19E-7   & 15.16   & 2.19E-7    & 15.30 \\
        $257$     & --        & --       & 1.51E-4  &  1.56E-10  & 14.96   & 1.56E-10   & 15.30 \\
          \hline
     \end{tabular}}                                                                                  
     \caption{$L^2$-errors when using the PCHIP, MQS, SPS, DBI, and PPI methods to approximate the function $f_{2}(x) = \frac{1}{1+ e^{-2kx}}, \quad k=100 \textrm{, and } x \in [-0.2,0.2]$.
           $N_{i}$ represents the number of input points used to build the approximation.
           $P_{j}$ represents the use of polynomials of degree $j$, with $j$ being the target degree for each interval.
           The sixth and eighth columns show the average polynomial degree used for the DBI and PPI methods respectively.
           The interval $[-0.2,0.2]$ is divided into $(N_{i}-1) / j$ elements and $j+1$ LGL quadrature points are used in each element.}
     \label{tab:heaviside_2}
   \end{table}

  \subsection{Example III $f_3(x)$}
  \label{subsec:example3r}
  The third example uses the modified function introduced in Equation \ref{eq:GT} with $\epsilon_{0}=0.01$, $\epsilon_{1}=1$ and $st=2$.
  The function $f_{3}(x)$ is particularly challenging because of the discontinuity at $x=-0.5$.
  The results from Tables \ref{tab:GelbT_1} and \ref{tab:GelbT_2} show that the $L^{2}$-errors from the four interpolation methods have the same order of accuracy.
  The DBI and PPI methods give slightly better approximation results compared to the other methods.
  The average polynomial degrees for the DBI and PPI approaches show that high-order polynomials are used.
  This suggest that in the smooth regions away from the discontinuity the DBI and PPI approaches lead to high-order accuracy.
  However, at the discontinuity, the DBI and PPI and other methods struggle to represent the underlying function.
  This example shows that both the DBI and PPI methods are appropriate approaches for interpolating from one mesh to another, 
because around the discontinuity, the methods are as accurate as the other ones, and in smooth regions the method gives better approximation results than the other approaches.
  The results from Tables \ref{tab:GelbT_1} and \ref{tab:GelbT_2} show that the $L^{2}$-error from the four interpolation methods have the same order.
    \begin{table}[H]
     \centering
     \scalebox{1.0}{
     \begin{tabular}{| c| c| c| c| c c| c c|}
          \hline
        $N_{i}$   & PCHIP    & MQS      & SPS   & \multicolumn{2}{c|}{DBI} & \multicolumn{2}{c|}{PPI}     \\
          \hline                              
                  & $L^{2}$-error & $L^{2}$-error & $L^{2}$-error & $L^{2}$-error & avg. deg. & $L^{2}$-error & avg. deg. \\
          \hline                              
                  &          &          &       & \multicolumn{4}{c|}{$\mathcal{P}_{1}$} \\
          \hline                          
        $17$      &  --      & --       & --    & 1.82E-1  & 1   & 1.82E-1 	&1  \\ 
        $33$      &  --      & --       & --    & 1.39E-1  & 1   & 1.39E-1 	&1  \\
        $65$      &  --      & --       & --    & 1.01E-1  & 1   & 1.01E-1 	&1  \\
        $129$     &  --      & --       & --    & 7.16E-2  & 1   & 7.16E-2 	&1  \\
        $257$     &  --      & --       & --    & 5.05E-2  & 1   & 5.05E-2 	&1  \\
          \hline                              
                  & $\mathcal{P}_{3}$  & $\mathcal{P}_{5}$ &  \multicolumn{5}{c|}{$\mathcal{P}_{4}$} \\
          \hline                          
        $17$      & 1.77E-1 & 1.59E-1  & 2.32E-1  & 1.82E-1   & 3.62   & 1.71E-1   & 3.81  \\ 
        $33$      & 1.39E-1 & 1.11E-1  & 1.56E-1  & 1.39E-1   & 3.97   & 1.29E-1   & 3.97  \\
        $65$      & 1.03E-1 & 7.90E-2  & 1.09E-1  & 9.38E-2   & 3.98   & 9.38E-2   & 3.98  \\
        $129$     & 7.42E-2 & 5.63E-2  & 7.69E-2  & 6.69E-2   & 3.99   & 6.70E-2   & 3.99  \\
        $257$     & 5.28E-2 & 4.04E-2  & 5.45E-2  & 4.73E-2   & 4      & 4.74E-2   & 4     \\
          \hline                                                                              
                  &          &          &       \multicolumn{5}{c|}{$\mathcal{P}_{8}$} \\
          \hline                                     
  	$17$      &--        &--        & 2.25E-1  & 1.83E-1   & 6.62   & 1.70E-1   & 7.06  \\
  	$33$      &--        &--        & 1.53E-1  & 1.36E-1   & 7.81   & 1.31E-1   & 7.81  \\
  	$65$      &--        &--        & 1.07E-1  & 9.62E-2   & 7.92   & 9.65E-2   & 7.92  \\
  	$129$     &--        &--        & 7.58E-2  & 6.90E-2   & 7.96   & 6.93E-2   & 7.96  \\
  	$257$     &--        &--        & 5.37E-2  & 4.90E-2   & 7.98   & 4.92E-2   & 7.98  \\
          \hline                                                                              
                  &          &          &        \multicolumn{5}{c|}{$\mathcal{P}_{16}$} \\
          \hline                                     
  	$17$      &--        &--        & 2.25E-1  & 1.82E-1   & 12.06   & 1.66E-1   & 13.06\\
  	$33$      &--        &--        & 1.50E-1  & 1.37E-1   & 14.09   & 1.33E-1   & 14.19 \\
  	$65$      &--        &--        & 1.05E-1  & 9.75E-2   & 15.78   & 9.81E-2   & 15.78 \\
  	$129$     &--        &--        & 7.45E-2  & 7.02E-2   & 15.90   & 7.07E-2   & 15.90 \\
  	$257$     &--        &--        & 5.29E-2  & 4.99E-2   & 15.95   & 5.03E-2   & 15.95 \\
          \hline
     \end{tabular}}                                                                                  
     \caption{$L^2$-errors when using the PCHIP, MQS, SPS, DBI, and PPI methods to approximate the function $f_{3}(x)$.
           $N_{i}$ represents the number of input points used to build the approximation.
           $P_{j}$ represents the use of polynomials of degree $j$, with $j$ being the target degree for each interval.
           The fifth and seventh columns show the average polynomial degree used for the DBI and PPI methods, respectively.
           The input points are uniformly distributed over the interval $[-1,1]$.} 
     \label{tab:GelbT_1}
   \end{table}
   \begin{table}[H]
     \centering
     \scalebox{1.0}{
     \begin{tabular}{| c| c| c| c| c c| c c|}
          \hline
        $N_{i}$   & PCHIP    & MQS      & SPS   &\multicolumn{2}{c|}{DBI} & \multicolumn{2}{c|}{PPI}     \\
          \hline                         
                  & $L^{2}$-error & $L^{2}$-error & $L^{2}$-error & $L^{2}$-error & avg. deg. & $L^{2}$-error & avg. deg. \\
          \hline                              
                  &          &          &       & \multicolumn{4}{c|}{$\mathcal{P}_{1}$} \\
          \hline                          
        $17$      &  --      & --       & --    &  1.82E-1   & 1   &  1.82E-1 	&1  \\ 
        $33$      &  --      & --       & --    &  1.39E-1   & 1   &  1.39E-1 	&1  \\
        $65$      &  --      & --       & --    &  1.01E-1   & 1   &  1.01E-1 	&1  \\
        $129$     &  --      & --       & --    &  7.16E-2   & 1   &  7.16E-2 	&1  \\
        $257$     &  --      & --       & --    &  5.05E-2   & 1   &  5.05E-2 	&1  \\
          \hline                              
                  & $\mathcal{P}_{3}$   & $\mathcal{P}_{5}$ &     \multicolumn{5}{c|}{$\mathcal{P}_{4}$} \\
          \hline                          
        $17$      & 1.64E-1 & 1.39E-1  & 1.87E-1  &  1.61E-1   & 3.81   & 1.58E-1   & 3.81   \\
        $33$      & 1.20E-1 & 9.79E-2  & 1.29E-1  &  1.18E-1   & 3.97   & 1.18E-1   & 3.97   \\
        $65$      & 8.70E-2 & 6.94E-2  & 9.04E-2  &  8.32E-2   & 3.98   & 8.53E-2   & 3.98   \\
        $129$     & 6.21E-2 & 4.96E-2  & 6.39E-2  &  5.93E-2   & 3.99   & 6.08E-2   & 3.99   \\
        $257$     & 4.39E-2 & 3.58E-2  & 4.54E-2  &  4.19E-2   & 4.00   & 4.29E-2   & 4      \\
          \hline                                                                              
                  &          &          &        \multicolumn{5}{c|}{$\mathcal{P}_{8}$} \\
          \hline                                     
        $17$      & --       & --       & 2.84E-1  &  1.85E-1   & 7.38   & 1.81E-01   & 7.50   \\
        $33$      & --       & --       & 9.61E-2  &  9.38E-2   & 7.62   & 1.27E-01   & 7.66   \\
        $65$      & --       & --       & 6.79E-2  &  6.75E-2   & 7.92   & 9.33E-02   & 7.92   \\
        $129$     & --       & --       & 4.82E-2  &  4.79E-2   & 7.96   & 6.72E-02   & 7.96   \\
        $257$     & --       & --       & 3.44E-2  &  3.37E-2   & 7.98   & 4.77E-02   & 7.98   \\
          \hline                                                                              
                  &          &          & \multicolumn{5}{c|}{$\mathcal{P}_{16}$} \\
          \hline                                     
        $17$     & --       & --       & 1.11E-1  &  1.08E-1   & 11.62   & 1.51E-1   & 12.12 \\
       $33$      & --       & --       & 1.86E-1  &  1.66E-1   & 14.31   & 1.55E-1   & 14.88 \\
       $65$      & --       & --       & 4.91E-2  &  5.06E-2   & 15.56   & 8.28E-2   & 15.58 \\
       $129$     & --       & --       & 3.51E-2  &  3.56E-2   & 15.90   & 5.92E-2   & 15.90 \\
       $257$     & --       & --       & 2.53E-2  &  2.48E-2   & 15.94   & 4.19E-2   & 15.94 \\
         \hline
    \end{tabular}}                                                                                  
    \caption{$L^2$-errors when using the PCHIP, MQS, SPS, DBI, and PPI methods to approximate the function $f_{3}(x)$.
          $N_{i}$ represents the number of input points used to build the approximation.
          $P_{j}$ represents the use of polynomials of degree $j$, with $j$ being the target degree for each interval.
          The sixth and eighth columns show the average polynomial degree used for the DBI and PPI methods respectively.
          The interval $[-1,1]$ is divided into $(N_{i}-1) / j$ elements and $j+1$ LGL quadrature points are used in each element.}
    \label{tab:GelbT_2}
  \end{table}

  \subsection{Example IV $f_4(x)$ }
  \label{subsec:example4r}
  The fourth example uses the function $f_{4}(x)$ defined in Equation \ref{eq:square} with $\epsilon_{0}=0.01$, $\epsilon_{1}=1$ and $st=2$.
  $f_{4}(x)$ depends on the size $h$ of each element, and as the number of element changes, so does the element size $h$ and the function $f_{4}(x)$.
  
  At the element boundaries $f_{4}(x)$ is only $C^{0}$-continuous with large gradients of opposite signs.
  The results from Tables \ref{tab:Square1} and \ref{tab:Square2} show that all the methods all the methods struggle to approximate the underlying function.
  With the exception of using a uniform mesh with PCHIP and SPS, the remaining results from Tables \ref{tab:Square1} and \ref{tab:Square2} show that all the methods have the same order of accuracy for both uniform and LGL meshes.
  The PPI and DBI methods give slightly smaller $L^{2}$-errors compared to the other approaches.
 
  \begin{table}[H]
     \centering
     \scalebox{1.0}{
     \begin{tabular}{| c| c| c| c| c| c c| c c|}
          \hline
        $N_{i}$   & PCHIP     & MQS      &$\mathcal{P}_{j}$ &  SPS & \multicolumn{2}{c|}{DBI} & \multicolumn{2}{c|}{PPI}     \\
          \hline                           
                  & $L^{2}$-error & $L^{2}$-error &    & $L^{2}$-error & $L^{2}$-error & avg. deg. & $L^{2}$-error & avg. deg. \\
          \hline                           
                  &    \multicolumn{8}{c|}{$ne = 4$} \\
          \hline                           
        $17$      &  3.74E-01 & 3.95E-1 & $\mathcal{P}_{4}$  &  3.72E-1 & 3.49E-1   & 2.62   & 3.49E-1   & 2.88   \\
        $33$      &  2.47E-01 & 2.59E-1 & $\mathcal{P}_{8}$  &  2.46E-1 & 2.19E-1   & 5.19   & 2.19E-1   & 6.19   \\
        $65$      &  1.55E-01 & 1.63E-1 & $\mathcal{P}_{16}$ &  1.54E-1 & 1.32E-1   & 13.34  & 1.32E-1   & 15.34  \\
          \hline                                                                                                  
                  &     \multicolumn{8}{c|}{$ne = 8$} \\                  
          \hline                                                          
        $33$      &  3.84E-01 & 3.94E-1 & $\mathcal{P}_{4}$  & 3.83E-1 & 4.04E-1   & 2.69   & 4.04E-1   & 2.94   \\
        $65$      &  2.54E-01 & 2.59E-1 & $\mathcal{P}_{8}$  & 2.52E-1 & 2.60E-1   & 5.34   & 2.74E-1   & 6.34   \\
        $129$     &  1.61E-01 & 8.23E-2 & $\mathcal{P}_{16}$ & 1.57E-1 & 1.67E-1   & 13.55  & 1.78E-1   & 15.30  \\
          \hline                                                                                                  
                  &     \multicolumn{8}{c|}{$ne = 16$} \\                 
          \hline                                                         
        $65$      & 3.90E-01  & 3.93E-1 & $\mathcal{P}_{4}$  &  3.89E-1 & 3.61E-1   & 2.72   & 3.61E-1   & 2.97   \\
        $129$     & 2.58E-01  & 2.58E-1 & $\mathcal{P}_{8}$  &  2.55E-1 & 2.26E-1   & 5.42   & 2.26E-1   & 6.42   \\
        $257$     & 1.63E-01  & 8.06E-2 & $\mathcal{P}_{16}$ &  1.58E-1 & 1.36E-1   & 13.65  & 1.36E-1   & 15.65  \\
          \hline              
     \end{tabular}}                                                                     
     \caption{$L^2$-errors when using the PCHIP, MQS, SPS, DBI, and PPI methods to approximate the Runge function $f_{4}(x)$.
           $N_{i}$ represents the number of input points used to build the approximation.
           $P_{j}$ represents the use of polynomials of degree $j$, with $j$ being the target degree for each interval.
           The value $ne$ represents the number of elements.
           The seventh and ninth columns show the average polynomial degree used for the DBI and PPI methods, respectively.
           The input points are uniformly distributed over the interval $[0,1]$.} 
     \label{tab:Square1}
   \end{table}
   \begin{table}[H]
     \centering
     \scalebox{1.0}{
     \begin{tabular}{| c| c| c| c| c| c c| c c|}
          \hline
        $N_{i}$   & PCHIP     & MQS    & $\mathcal{P}_{j}$ & SPS & \multicolumn{2}{c|}{DBI} & \multicolumn{2}{c|}{PPI}     \\
          \hline                           
                  & $L^{2}$-error & $L^{2}$-error &    & $L^{2}$-error & $L^{2}$-error & avg. deg. & $L^{2}$-error & avg. deg. \\
          \hline                           
                  &    \multicolumn{8}{c|}{$ne = 4$} \\
          \hline                           
        $17$      & 3.02E-1  & 3.18E-1 & $\mathcal{P}_{4}$  & 3.02E-1 & 2.32E-1   & 2.75   & 2.32E-1   & 3.00   \\
        $33$      & 1.33E-1  & 1.39E-1 & $\mathcal{P}_{8}$  & 1.33E-1 & 9.60E-2   & 5.25   & 9.60E-2   & 6.25   \\
        $65$      & 3.80E-2  & 3.92E-2 & $\mathcal{P}_{16}$ & 3.77E-2 & 2.11E-2   & 11.88  & 2.11E-2   & 13.31  \\
          \hline                                                                                                                        
                  &    \multicolumn{8}{c|}{$ne = 8$} \\                                         
          \hline                                                                                
        $33$      & 3.10E-1  & 3.17E-1 & $\mathcal{P}_{4}$  & 3.10E-1 & 3.42E-1   & 2.75   & 3.42E-1   & 3.00   \\
        $65$      & 1.37E-1  & 1.39E-1 & $\mathcal{P}_{8}$  & 1.36E-1 & 1.59E-1   & 5.25   & 1.65E-1   & 6.12   \\
        $129$     & 3.98E-2  & 3.72E-2 & $\mathcal{P}_{16}$ & 3.87E-2 & 5.33E-2   & 11.88  & 5.60E-2   & 13.31  \\
          \hline                                                                                                                           
                  &    \multicolumn{8}{c|}{$ne = 16$} \\                                        
          \hline                                                                                
        $65$      & 3.14E-1  & 3.16E-1 & $\mathcal{P}_{4}$  & 3.14E-1 & 2.32E-1   & 2.75   & 2.32E-1   & 3.00  \\
        $129$     & 1.39E-1  & 1.38E-1 & $\mathcal{P}_{8}$  & 1.37E-1 & 9.60E-2   & 5.25   & 9.60E-2   & 6.25  \\
        $257$     & 4.07E-2  & 3.37E-2 & $\mathcal{P}_{16}$ & 3.90E-2 & 2.11E-2   & 11.88  & 2.11E-2   & 13.31 \\
          \hline              
     \end{tabular}}                                                                     
     \caption{$L^2$-errors when using the PCHIP, MQS, SPS, DBI, and PPI methods to approximate the Runge function $f_{4}(x)$.
           $N_{i}$ represents the number of input points used to build the approximation.
           $P_{j}$ represents the use of polynomials of degree $j$, with $j$ being the target degree for each interval.
           The value $ne$ represents the number of elements.
           The seventh and ninth columns show the average polynomial degree used for the DBI and PPI methods respectively.
           The interval $[0,1]$ is divided into $(N_{i}-1) / j$ elements and $j+1$ LGL quadrature points are used in each element.}
     \label{tab:Square2}
   \end{table}

  \subsection{Example V $f_5(x)$}
  \label{subsec:example5r}
  The fifth experiment uses the function $f_{5}(x)$ defined in Equation \ref{eq:tanh} with $\epsilon_{0}=0.01$, $\epsilon_{1}=1$ and $st=1$.
  $f_{5}(x)$ depends on the size $h$ of each element and as the number of elements changes, so does the element size $h$ and $f_{5}(x)$.
  Similarly to $f_{4}(x)$, $f_{5}(x)$ is only $C^{0}$-continuous at the element boundaries.
  However, the gradients remain positive over the entire interval.
  This example is constructed to reflect the $C^{0}$-continuity observed in the spectral element method used in NEPTUNE.
  Tables \ref{tab:Tanh1} and \ref{tab:Tanh2} shows that the approximation errors from PCHIP, MQS, and SPS methods improve slowly compared to the DBI and PPI methods, as we increase the polynomial degree and the number of points.
  The PCHIP, SPS, and MQS methods use approximations of the first derivatives and  enforce $C^{1}$-continuity at the element boundaries. 
  Overall, the results from Tables \ref{tab:Tanh1} and \ref{tab:Tanh2} show that the DBI and PPI methods has smaller $L^{2}$-errors compared to the remaining methods. 
  \begin{table}[H]
     \centering
     \scalebox{1.0}{
     \begin{tabular}{| c| c| c| c| c| c c| c c|}
          \hline
        $N_{i}$   & PCHIP     & MQS      & $\mathcal{P}_{j}$ & SPS & \multicolumn{2}{c|}{DBI} & \multicolumn{2}{c|}{PPI}     \\
          \hline                           
                  & $L^{2}$-error & $L^{2}$-error &    & $L^{2}$-error & $L^{2}$-error & avg. deg. & $L^{2}$-error & avg. deg. \\
          \hline                           
                  &    \multicolumn{8}{c|}{$ne = 4$} \\
          \hline                           
        $17$      & 1.68E-2  & 4.50E-2  & $\mathcal{P}_{4}$  &  1.23E-2 & 2.36E-2   & 3.56   & 2.36E-2   & 3.56  \\
        $33$      & 9.95E-3  & 6.04E-3  & $\mathcal{P}_{8}$  &  1.36E-2 & 4.20E-4   & 7.59   & 4.20E-4   & 7.62  \\
        $65$      & 1.67E-3  & 3.82E-4  & $\mathcal{P}_{16}$ &  5.77E-3 & 3.65E-5   & 14.98  & 3.65E-5   & 14.98 \\
          \hline                                                                                                  
                  &    \multicolumn{8}{c|}{$ne = 8$} \\                                      
          \hline                                                                             
        $33$      & 1.99E-2  & 1.20E-2  & $\mathcal{P}_{4}$  &  1.29E-2 & 1.65E-2   & 3.91   & 1.65E-2   & 3.91  \\
        $65$      & 3.35E-3  & 7.63E-4  & $\mathcal{P}_{8}$  &  7.42E-3 & 2.17E-4   & 7.67   & 2.17E-4   & 7.67  \\
        $129$     & 3.70E-4  & 4.44E-5  & $\mathcal{P}_{16}$ &  2.89E-3 & 5.01E-5   & 14.84  & 5.01E-5   & 14.88 \\
          \hline                                                                                                  
                  &    \multicolumn{8}{c|}{$ne = 16$} \\                                     
          \hline                                                                             
        $33$      & 6.73E-3 & 2.46E-3   & $\mathcal{P}_{4}$  &  4.57E-3 & 1.18E-3   & 3.86   & 1.18E-3   & 3.86  \\
        $65$      & 8.22E-4 & 4.53E-4   & $\mathcal{P}_{8}$  &  3.61E-3 & 5.27E-5   & 7.51   & 5.27E-5   & 7.51  \\
        $256$     & 1.53e-4 & 1.59E-4   & $\mathcal{P}_{16}$ &  1.42E-3 & 4.27E-11  & 14.65  & 4.27E-11   & 14.70 \\
          \hline              
     \end{tabular}}                                                                     
     \caption{$L^2$-errors when using the PCHIP, MQS, SPS, DBI, and PPI methods to approximate the function $f_{5}(x)$.
           $N_{i}$ represents the number of input points used to build the approximation.
           $P_{j}$ represents the use of polynomials of degree $j$, with $j$ being the target degree for each interval.
           The value $ne$ is the number of elements used.
           The seventh and ninth columns show the average polynomial degree used for the DBI and PPI methods respectively.
           The input points are uniformly distributed over the interval $[-2,0]$.} 
     \label{tab:Tanh1}
  \end{table}
  \begin{table}[H]
     \centering
     \scalebox{1.0}{
     \begin{tabular}{| c| c| c| c| c| c c| c c|}
          \hline
        $N_{i}$   & PCHIP     & MQS    & $\mathcal{P}_{j}$  & SPS & \multicolumn{2}{c|}{DBI} & \multicolumn{2}{c|}{PPI}     \\
          \hline                           
                  & $L^{2}$-error & $L^{2}$-error &    & $L^{2}$-error & $L^{2}$-error & avg. deg. & $L^{2}$-error & avg. deg. \\
          \hline                           
                  &    \multicolumn{8}{c|}{$ne = 4$} \\
          \hline                           
        $17$      & 4.64E-2 & 3.10E-2  & $\mathcal{P}_{4}$  & 1.23E-2 & 5.09E-2   & 3.06   & 5.09E-2   & 3.06  \\
        $33$      & 7.43E-3 & 4.29E-3  & $\mathcal{P}_{8}$  & 1.36E-2 & 1.47E-3   & 6.34   & 1.47E-3   & 6.44  \\
        $65$      & 9.48E-4 & 1.58E-4  & $\mathcal{P}_{16}$ & 5.77E-3 & 3.15E-6   & 14.83  & 3.15E-6   & 14.84  \\
          \hline                                                                                                
                  &    \multicolumn{8}{c|}{$ne = 8$} \\                                      
          \hline                                                                             
        $33$      & 2.57E-2 & 9.45E-3  & $\mathcal{P}_{4}$  & 1.29E-2 & 1.52E-2   & 3.84   & 1.52E-2   & 3.84 \\
        $65$      & 3.18E-3 & 9.15E-4  & $\mathcal{P}_{8}$  & 7.42E-3 & 2.66E-4   & 7.66   & 2.66E-4   & 7.67 \\
        $129$     & 4.08E-4 & 2.71E-5  & $\mathcal{P}_{16}$ & 2.89E-3 & 3.75E-4   & 14.88  & 4.53E-4   & 14.88 \\
          \hline                                                                                                
                  &    \multicolumn{8}{c|}{$ne = 16$} \\                                     
          \hline                                                                             
        $33$      & 8.03E-3 & 3.77E-3  & $\mathcal{P}_{4}$  &  4.57E-3 & 1.91E-3   & 3.81   & 1.91E-3   & 3.81  \\
        $65$      & 9.93E-4 & 2.22E-4  & $\mathcal{P}_{8}$  &  3.61E-3 & 2.23E-6   & 7.51   & 2.23E-6   & 7.51  \\
        $256$     & 1.23E-4 & 3.30E-5  & $\mathcal{P}_{16}$ &  1.42E-3 & 2.63E-12  & 14.37  & 2.63E-12  & 14.50  \\
          \hline              
     \end{tabular}}                                                                     
     \caption{$L^2$-errors when using the PCHIP, MQS, SPS, DBI, and PPI methods to approximate the function $f_{5}(x)$.
           $N_{i}$ represents the number of input points used to build the approximation.
           $P_{j}$ represents the use of polynomials of degree $j$, with $j$ being the target degree for each interval.
           The value of $ne$ is the number of elements used.
           The seventh and ninth columns show the average polynomial degree used for the DBI and PPI methods respectively.
           The interval $[-2,0]$ is divided into $(N_{i}-1) / j$ elements and $j+1$ LGL quadrature points are used in each element.}
     \label{tab:Tanh2}
   \end{table}

  \subsection{Example VI BOMEX}
  \label{subsec:bomex}
      %
      The 1D Barbados Oceanographic and Meteorological Experiment (BOMEX) \cite{friedman1970essa} is a single column test case that was developed to measure and study changes in the properties of heat, moisture, and momentum.
      In this example, the dynamics and physics results are calculated on different meshes.
      The dynamics uses uniformly spaced points that indicate the boundary of each level in the vertical column.
      The physics mesh is constructed using the mid-point of each level.
      The advections in the dynamics are approximated using a fifth order weighted essentially non-oscillatory (WENO) and third-order Runge-Kutta methods ~\cite{doi:10.1080/1061856031000104851}.
      At each time step, the dynamics are calculated on the dynamics mesh, and the results are interpolated to the physics mesh for the use of the physics routines.
      The physics terms are calculated using the physics mesh, and the results are interpolated back to the dynamics mesh.

\begin{figure}[H]
         \centering
         \includegraphics[scale = 0.15]{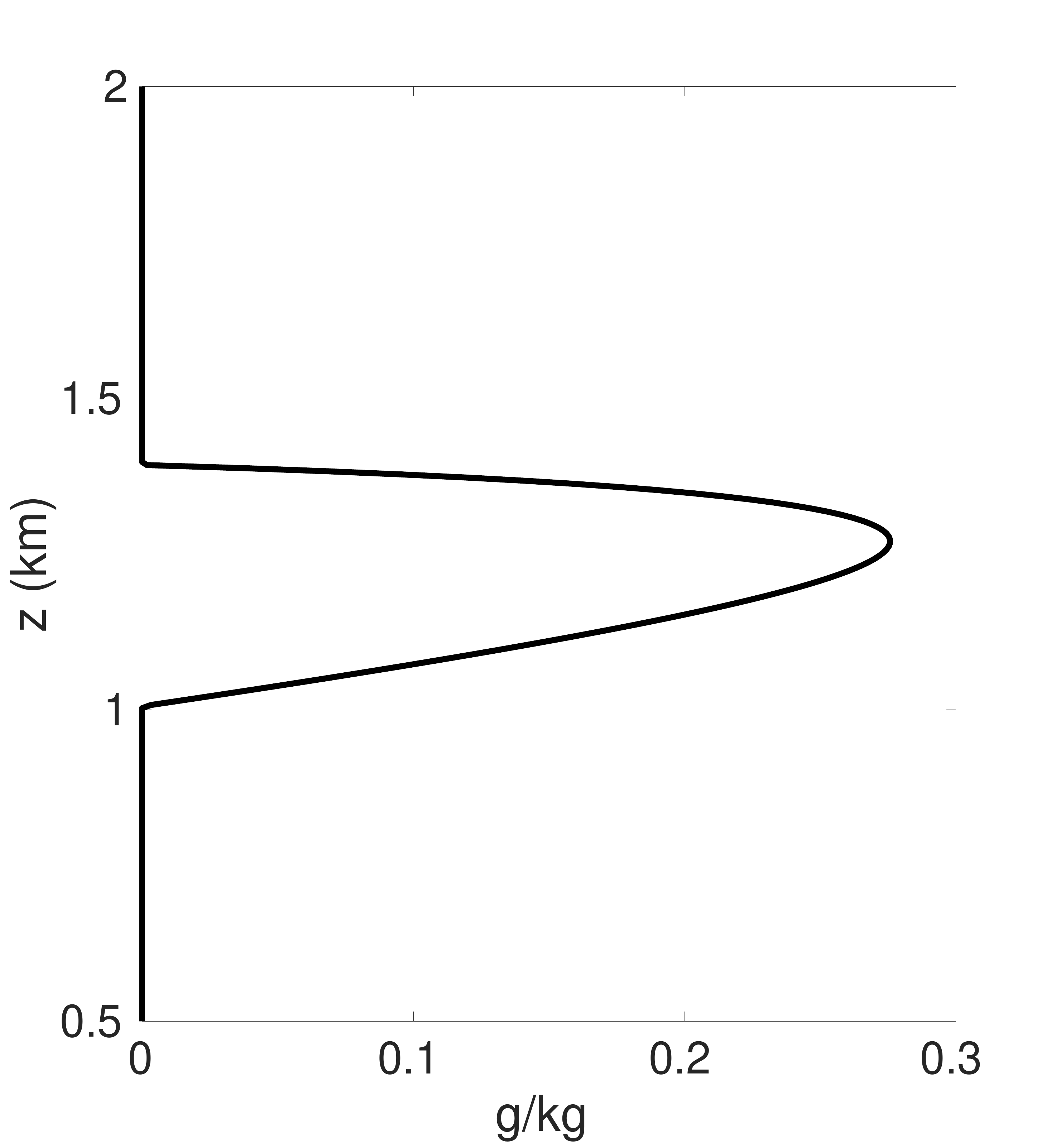}
         \caption{Target cloud mixing ratio $q_{c}$ profile  from  BOMEX test case at $t=5 h$ with $nz = 600$ points.
                  A fifth-order WENO  and third-order Runge-Kutta schemes with $CFL= 0.1$ are used for the dynamics (advection).
                  The same mesh is used for the dynamics and physics calculations.}
         \label{fig:qcweno_physics}
      \end{figure}

\begin{figure}[H]
         \centering
         \includegraphics[scale = 0.15]{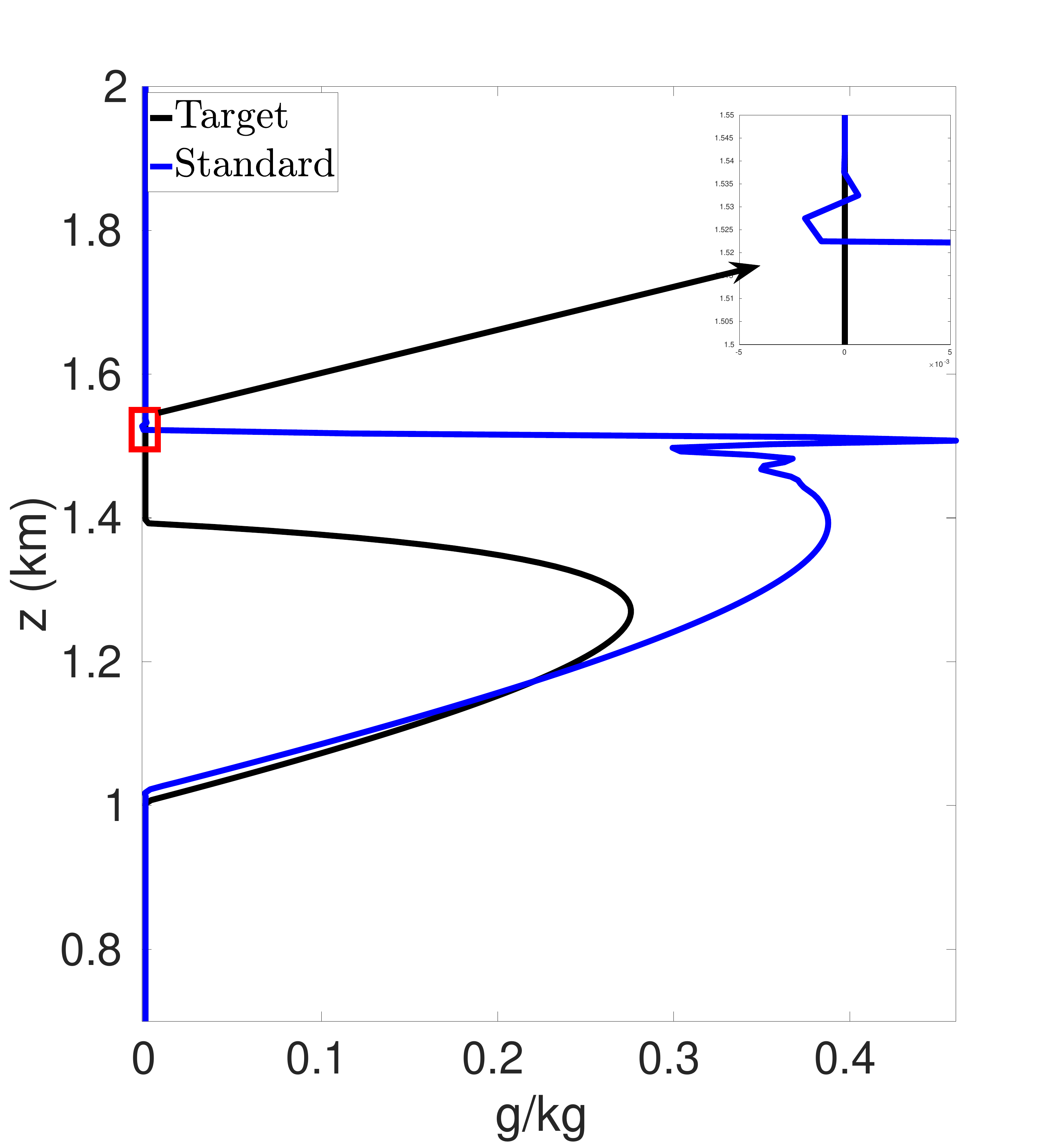}
         \caption{Cloud mixing ratio $q_{c}$ profile  from  BOMEX test case at $t= 5 h$ with $nz=600$ points. 
                  The profile in black is the target solution. 
                  The profile in blue is obtained using a standard interpolation method when mapping solution values between dynamics and physics meshes.
                  A fifth-order WENO  and third-order Runge-Kutta schemes with $CFL= 0.1$ are used for the dynamics (advection).}
         \label{fig:qcweno_std}
      \end{figure}

\begin{figure}[H]
         \centering
         \includegraphics[scale = 0.15]{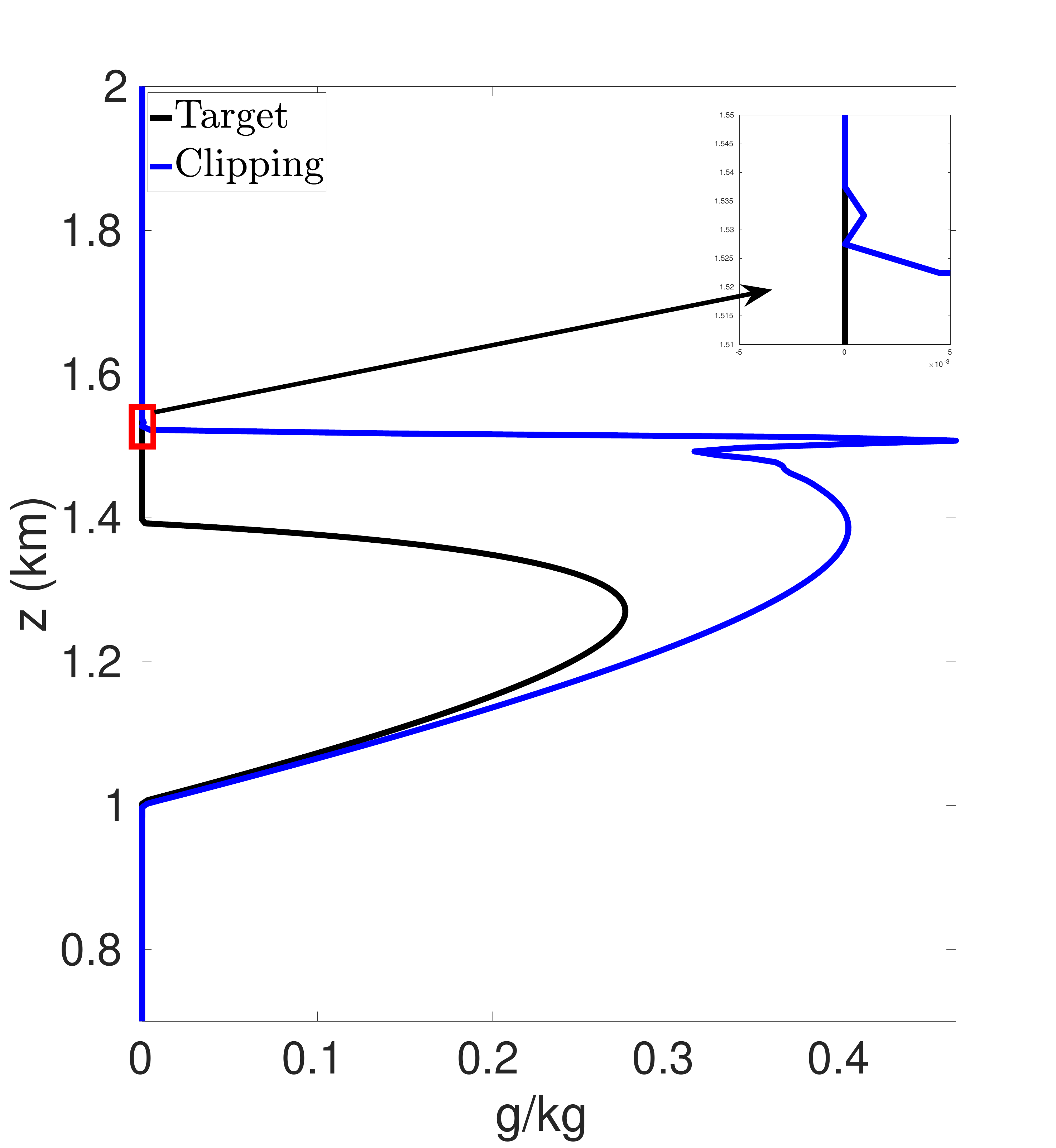}
         \caption{Cloud mixing ratio $q_{c}$ profile from  BOMEX test case at $t = 5 h$ with $nz=601$ points. 
                  The profile in black is the target solution. 
                  The profile in blue is obtained using a standard interpolation method when mapping solution values between dynamics and physics meshes.
                  ``Clipping is used after the interpolation to remove any negative value.
                  A fifth-order WENO  and third-order Runge-Kutta schemes with $CFL= 0.1$ are used for the dynamics (advection).}
         \label{fig:qcweno_clip}
      \end{figure}

\begin{figure}[H]
         \centering
         \includegraphics[scale = 0.15]{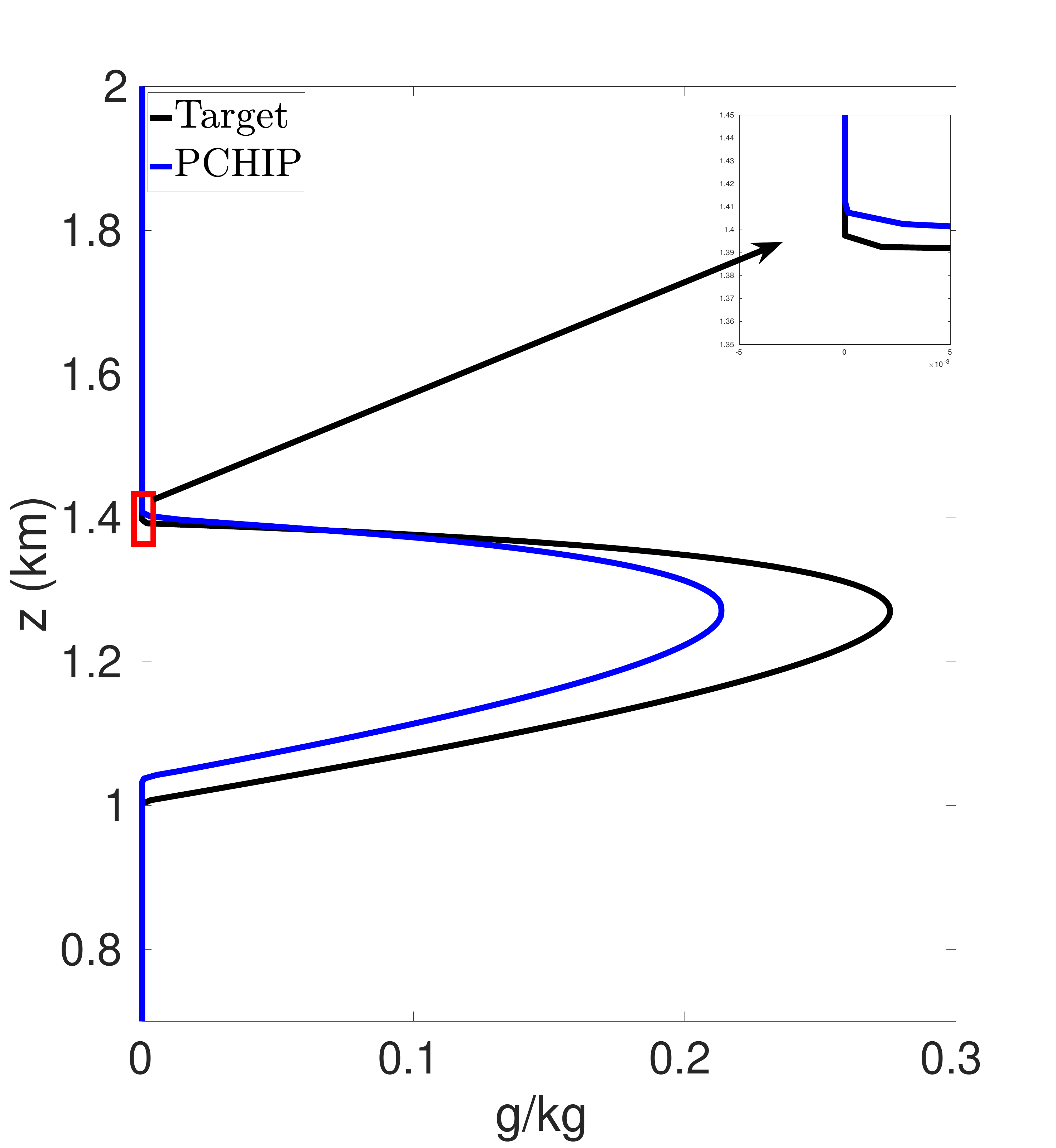}
         \caption{Cloud mixing ratio $q_{c}$ profile  from  BOMEX test case at $t= 5 h$ with $nz= 600$ points. 
                  The profile in black is the target solution. 
                  The profile in blue is obtained using the PCHIP method when mapping solution values between dynamics and physics meshes.
                  A fifth-order WENO  and third-order Runge-Kutta schemes with $CFL= 0.1$ are used for the dynamics (advection).}
         \label{fig:qcweno_pchip}
      \end{figure}

\begin{figure}[H]
         \begin{subfigure}{0.5\textwidth}
         \centering
         \includegraphics[scale = 0.11]{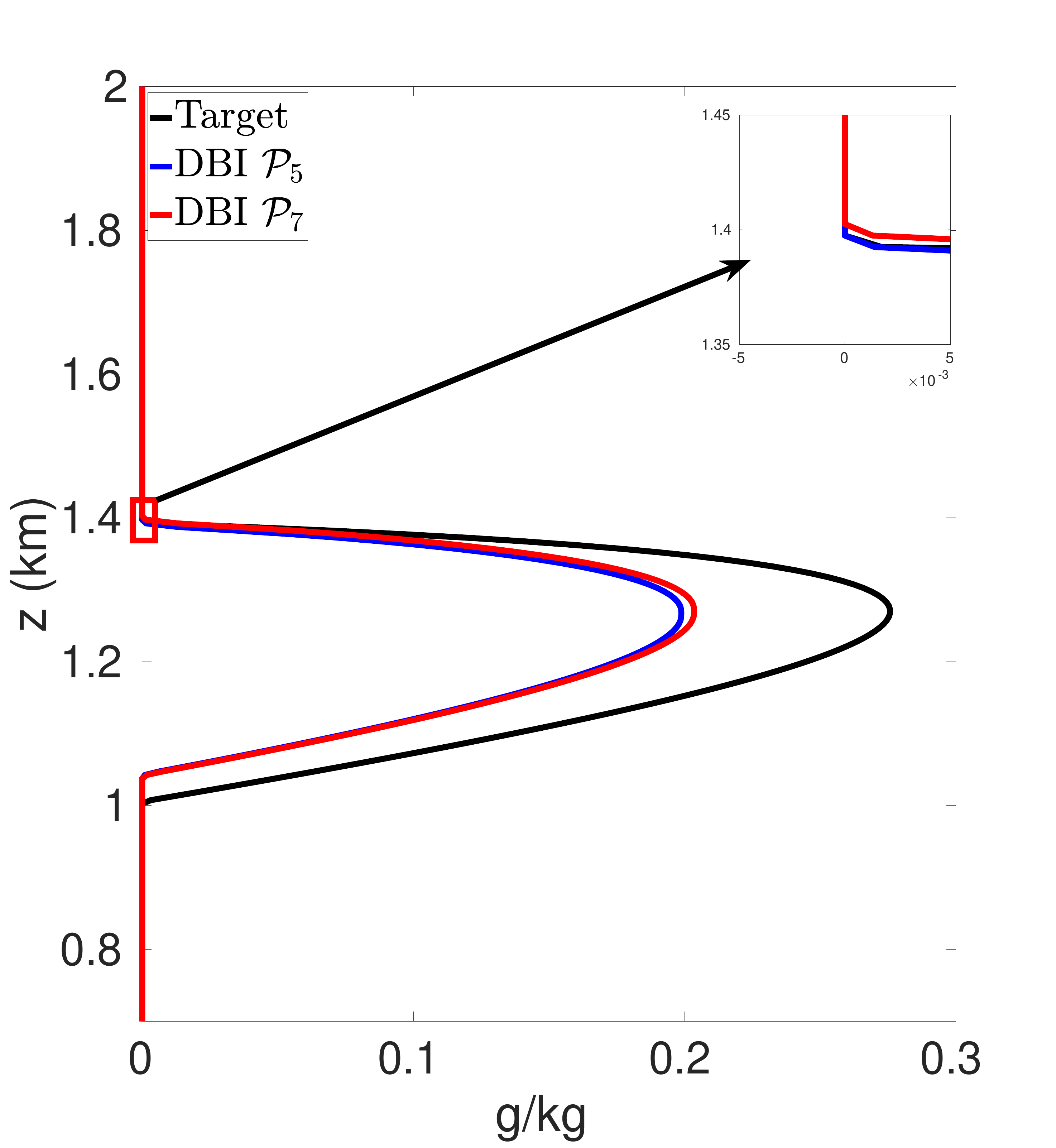}
         \subcaption{}
         \end{subfigure}
         \begin{subfigure}{0.5\textwidth}
         \centering
         \includegraphics[scale = 0.11]{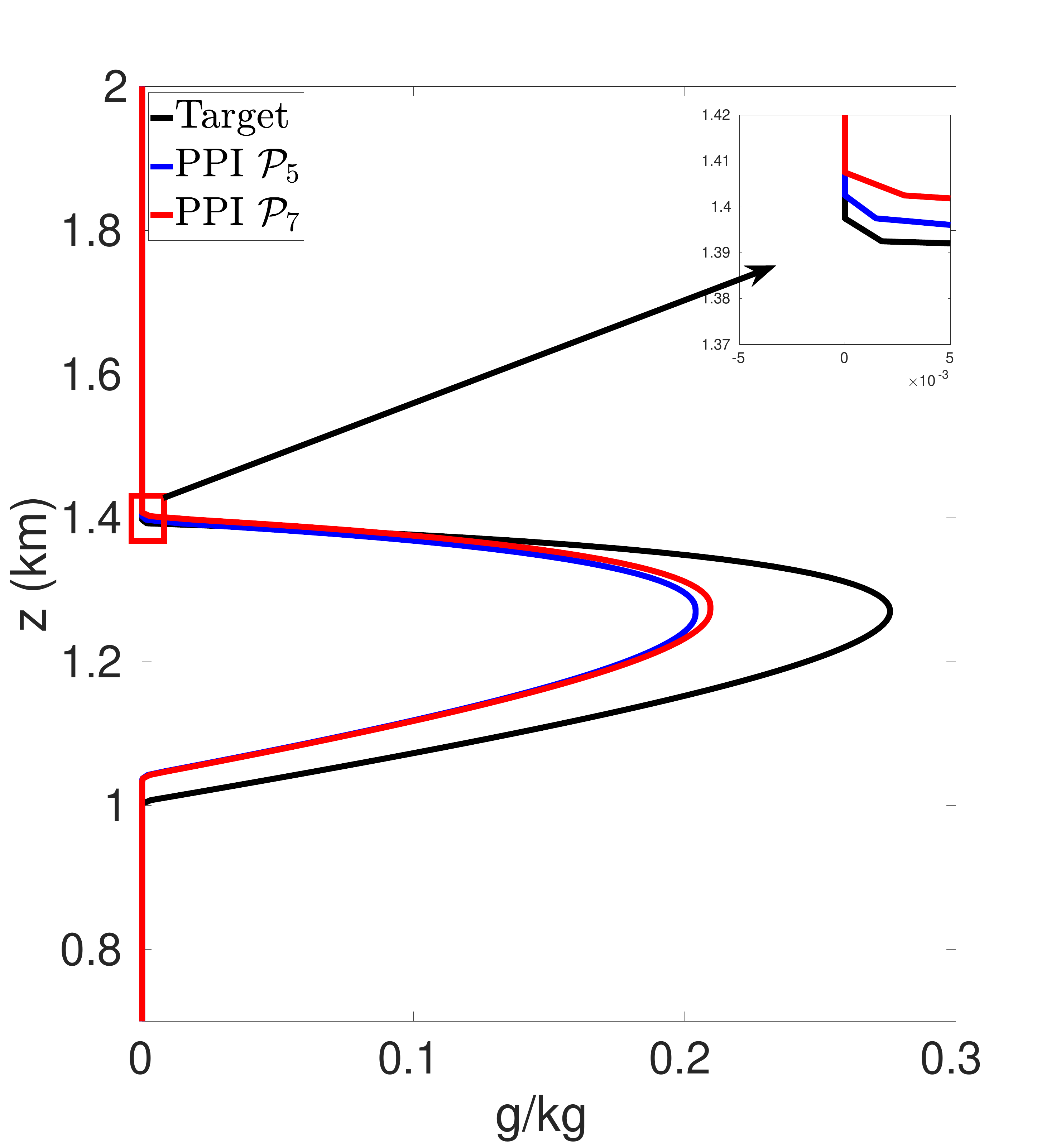}
         \subcaption{}
         \end{subfigure}
         \begin{subfigure}{0.5\textwidth}
         \centering
         \includegraphics[scale = 0.11]{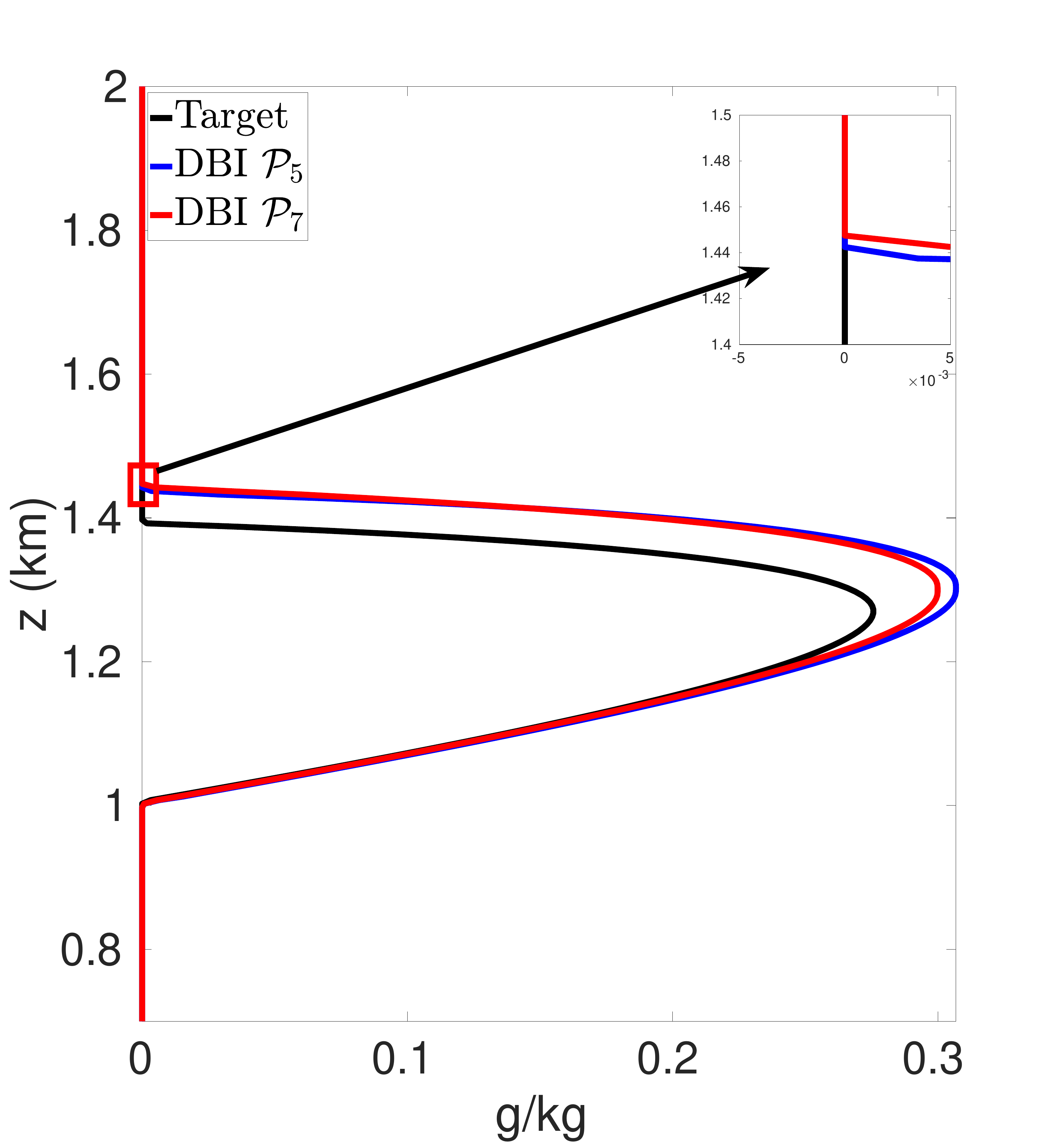}
         \subcaption{}
         \end{subfigure}
         \begin{subfigure}{0.5\textwidth}
         \centering
         \includegraphics[scale = 0.11]{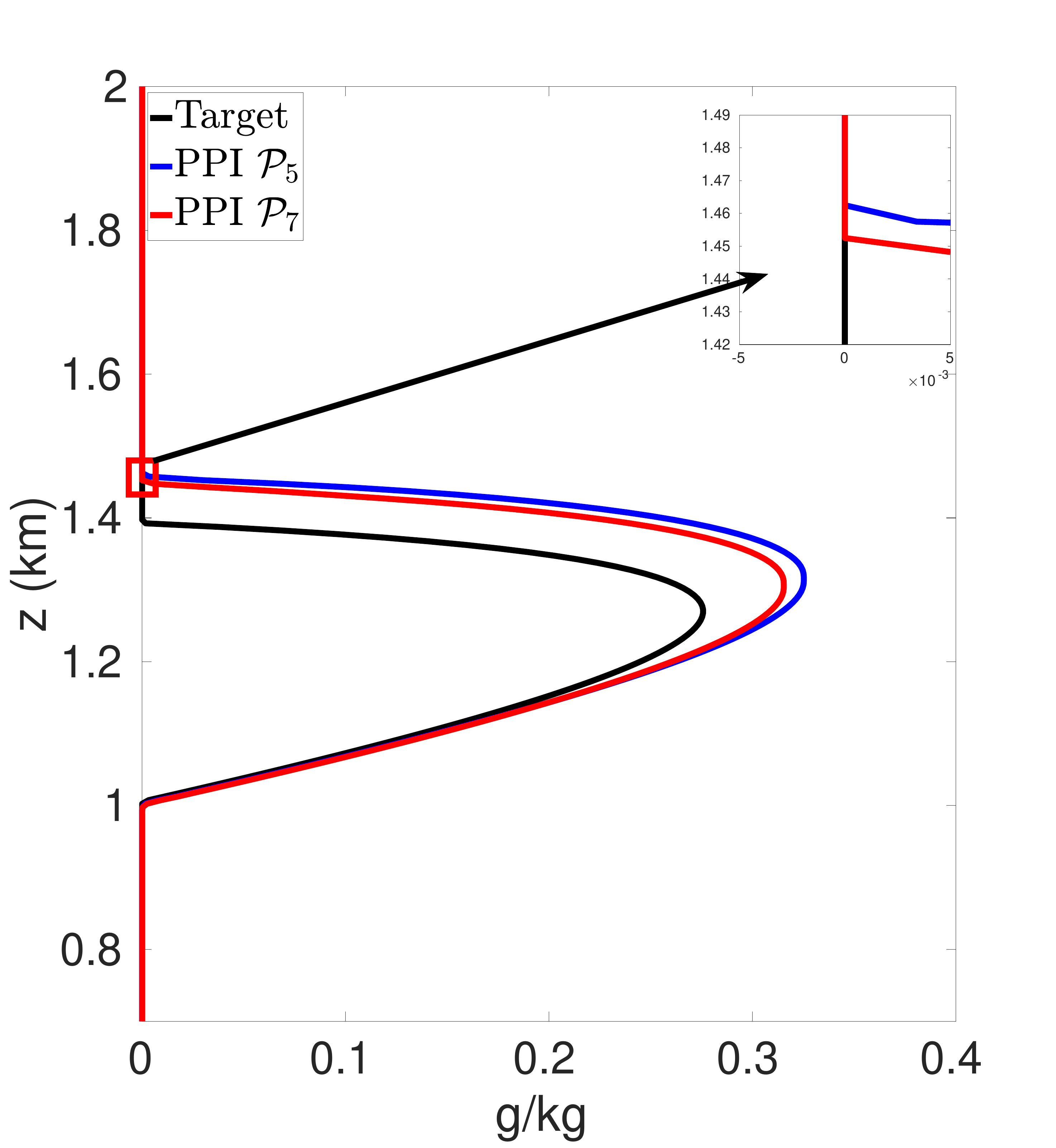}
         \subcaption{}
         \end{subfigure}
         \begin{subfigure}{0.5\textwidth}
         \centering
         \includegraphics[scale = 0.11]{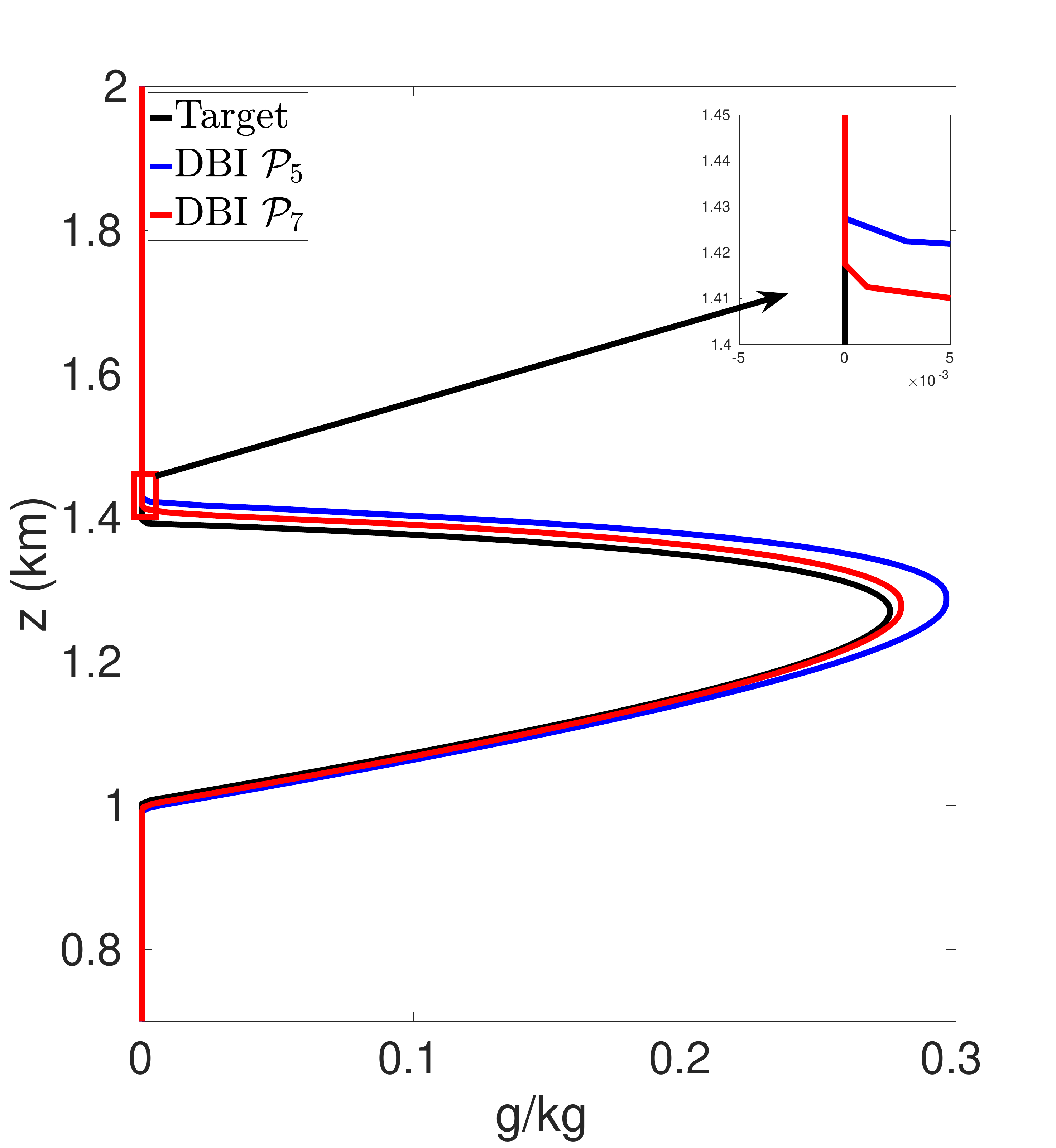}
         \subcaption{}
         \end{subfigure}
         \begin{subfigure}{0.5\textwidth}
         \centering
         \includegraphics[scale = 0.11]{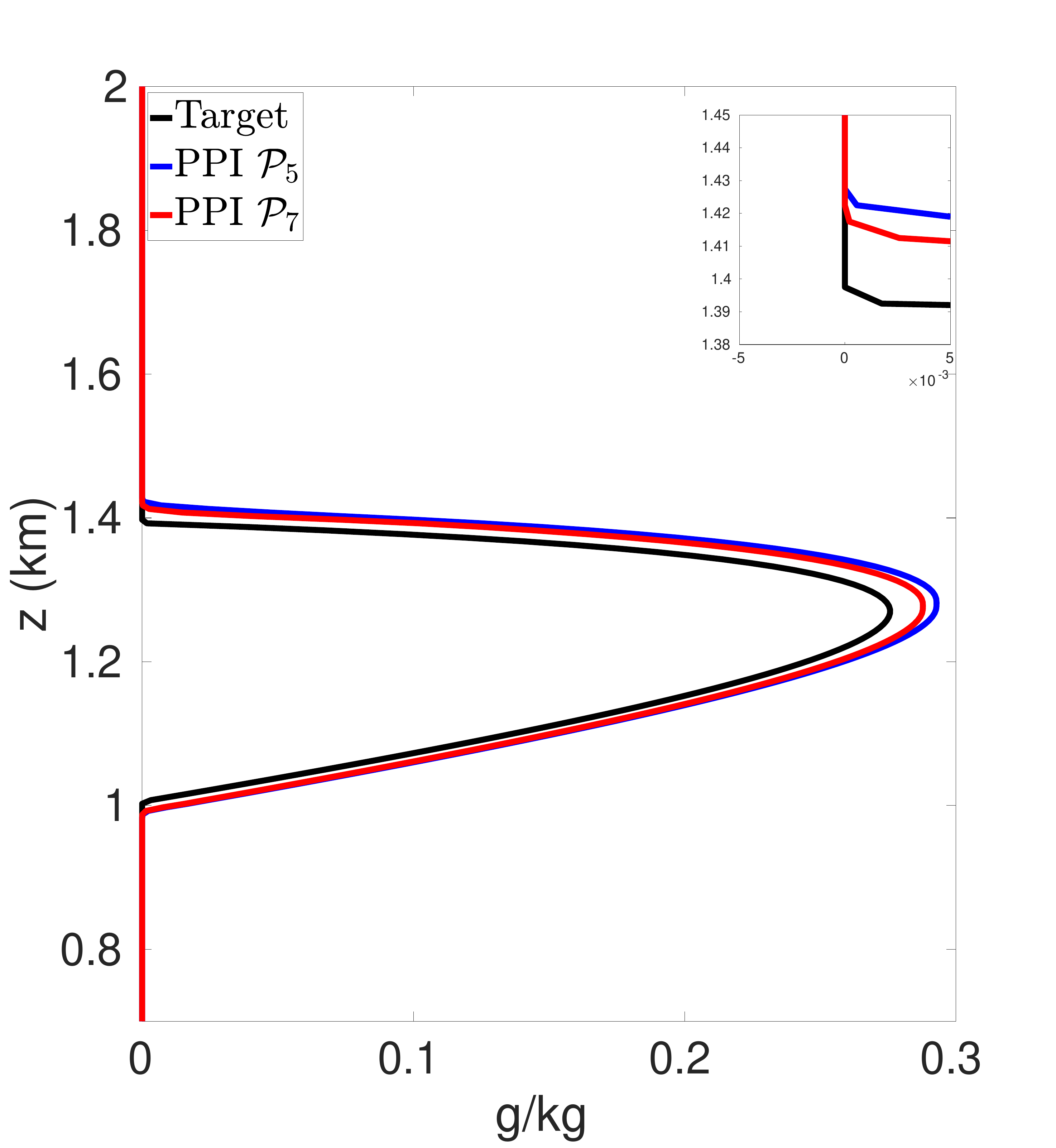}
         \subcaption{}
         \end{subfigure}
         \caption{Cloud mixing ratio $q_{c}$ profile  from  BOMEX test case at $t= 5 h$ with $nz= 600$ points with $\epsilon_{0}=\epsilon_{1}=10^{-5}$. 
                  The profile in black is the target solution. 
                  The profiles on the left and right are obtained using the DBI and PPI methods, respectively, to map solution values between meshes.
                  The maximum polynomial degrees are set to $5$ and $7$ for the blue and red plots, respectively.
                  A fifth-order WENO  and third-order Runge-Kutta schemes with $CFL= 0.1$ are used for the dynamics (advection).}
         \label{fig:qcweno_ppi}
      \end{figure}

      As in \cite{Rotstayn2000}, let $q_{c}$ be the cloud water mixing ratio profile in the different experiments.
      Figures \ref{fig:qcweno_physics}-\ref{fig:qcweno_ppi} show the cloud mixing ratio profile  $q_{c}$ at $t=5h$ that is used as input for the physics routines.
      The physics calculations require positive input values for $q_{c}$.
      Figure \ref{fig:qcweno_physics} shows the target profile for $q_{c}$.
      This target profile is obtained by using the same mesh for both dynamics and physics calculations where mapping is not required and $q_{c}$ remains positive during the simulation. 
      In addition, as the temporal and spatial resolution increases, $q_{c}$ converges to the profile shown in Figure \ref{fig:qcweno_physics}.
      Figures \ref{fig:qcweno_std}-\ref{fig:qcweno_ppi} are used to investigate different interpolation methods for mapping the solution values between meshes in the case where the dynamics and physics are calculated using different meshes. 

      Figure \ref{fig:qcweno_std} shows the cloud mixing ratio profiles $q_{c}$ for the target and approximated solution at $t= 5 h$.
      In the case of the approximated solution, a fifth-order standard polynomial interpolation is used when mapping between dynamics and physics meshes.
      For a given interval $I_{i}$, the polynomial interpolant is constructed using the stencil $\mathcal{V}_{4}=\{x_{i-2}, x_{i-1}, x_{i}, x_{i+1}, x_{i+2}, x_{i+3}, \}$.
      At the boundary and nearby boundary intervals, the stencil $\mathcal{V}_{4}$ is biased toward the interior of the domain. 
      The results in Figure \ref{fig:qcweno_std} demonstrate that using the standard polynomial interpolation lead to oscillations, negatives values, and an overestimation of the peak and total cloud mixing ratio of the profile $q_{c}$.
      Using standard polynomial interpolation leads to an overproduction of the total cloud mixing ratio by $93.45\%$.
      The peak is $max(q_{c}) = 0.46 g/kg$, which is larger than the target peak $max(q_{c})= 0.28 g/kg$. 

      The negative values in Figure \ref{fig:qcweno_std} can be removed via ``clipping".
      ``Clipping" is a procedure that consists of removing the negative values by setting them to zero.
      Figure \ref{fig:qcweno_clip} shows the cloud mixing ratio profiles for the target solution and an approximated solution that uses ``clipping" to remove the negative values at each time step.
      The approximated solution uses a standard interpolation to map the data values from one mesh to another. 
      The interpolant for each interval is constructed using the stencil $\mathcal{V}_{4}=\{x_{i-2}, x_{i-1}, x_{i}, x_{i+1}, x_{i+2}, x_{i+3}, \}$ with a fifth order polynomial.
      Once the interpolation is completed, ``clipping" is used to remove the negative values.
      Figure \ref{fig:qcweno_clip} shows that using ``clipping"  still allows for oscillations and a positive bias in the prediction of cloud mixing ratio $q_{c}$. 
      The total cloud mixing ratio is $2.09$ greater than the target solution and the peak $max(q_{c})=0.46 g/kg$ is larger than the target peak $max(q_{c})=0.28 g/kg$.

      Using PCHIP for the mapping between the dynamics and physics meshes eliminates the negative values, remove oscillations, and reduces the positive bias in the cloud mixing ratio prediction compared to the standard interpolation with and without ``clipping".
      Figure \ref{fig:qcweno_pchip} shows the target profile $q_{c}$ and an approximated profile that uses PCHIP for mapping solution values between dynamics and physics meshes.
      %
      %
      The total cloud mixing ratio is now $27.21\%$ less than the target with a peak $max(q_{c})= 0.21 g/kg$.
      In the case of the BOMEX test case, NEPTUNE, and similar codes, using PCHIP for mapping data values from one mesh to another can degrade the high-order accuracy obtained from the high-order methods used for the dynamics calculations.
      PCHIP is only third-order whereas the dynamics calculations use a fifth order method.
      This limitation can be addressed via high-order data-bounded or positivity-preserving methods. 

      Figure \ref{fig:qcweno_ppi} shows cloud mixing ratio profiles for the target and approximated solutions that use the DBI and PPI methods to map the solution values between meshes.
      The maximum polynomial degree for the DBI and PPI methods is set to $5$ and $7$, and the parameters $\epsilon_{0}$ and $\epsilon_{1}$ are both set a value of $10^{-5}$.
      For larger values of $\epsilon_{0}$ and $\epsilon_{1}$, the PPI approach introduces oscillations that lead to positive bias prediction of the cloud mixing ratio.
      These oscillation are caused by the relaxed nature of the PPI approach, which still allows the interpolants to oscillate while remaining positive.
      The positive bias and oscillations can be removed using the DBI or PPI method with small values for $\epsilon_{0}$ and $\epsilon_{1}$.
      When using the PPI method for mapping, the total amount of the cloud mixing ratio is less than the target for $st=1$ and more than the target for $st=2$ and $st=3$.
      The parameter $st$ is described in Section \ref{subsec:dbippi}.
      Figure \ref{fig:qcweno_ppi} shows that using the  DBI and PPI methods with $\epsilon_{0}=\epsilon_{1}=10^{-5}$ to map data values between the dynamics and physics meshes eliminates the negative values, removes the oscillations, and significantly reduces the positive bias in the cloud mixing ratio prediction.
      Using the DBI and PPI methods leads to better approximation of the peak value of the total cloud mixing ratio compared to using the standard interpolation and  PCHIP approaches. 
      The best approximation of the total amount of the cloud mixing ratio is with the DBI method, which is $7.57\%$ more than the target with a peak of $max(q_{c}) = 0.28  g/kg$.

      In summary, using DBI and PPI methods to map data values between both dynamics and physics meshes produces better approximation results compared to the standard interpolation and  PCHIP methods.
      Tables \ref{tab:qc0} and \ref{tab:qc1} provide a summary of the maximum values and the total amount of the cloud mixing ratios for each case.
      The DBI and PPI methods with a target polynomial set to $d=7$ lead to a better approximation of the peak and the total the cloud mixing ratios compared to the standard interpolation and PCHIP approaches.
      The results from Tables \ref{tab:qc0} and \ref{tab:qc1} indicate that the DBI method is the most suitable approach to map data values between meshes for the BOMEX test case.
      This study provided an example demonstrating how to use the DBI and PPI methods for mapping data values between meshes in the context of NWP.
      The BOMEX example also demonstrated that positivity alone may not be sufficient to remove oscillations in the solution, and the interpolants may need to be constrained to be between the data values for better approximation.
    \begin{table}[H]
      \centering
      \begin{tabular}{ c| c c c c}
        \hline
        \hline
                        & Target    & STD     & Clipping & PCHIP   \\
        \hline
        maximum $q_{c}$ & 0.28      & 0.46    & 0.46     & 0.21    \\
        total $q_{c}$   & 69.82   & 135.07   & 145.89   & 50.82    \\
        \hline
        \hline
      \end{tabular}                                                                                  
      \caption{Maximum values of $q_{c}$ and the total amount of the cloud mixing ratio at $t = 5 h$ with $nz=600$ points.
               The total amount of the cloud mixing ratio is calculated by estimating the integral $q_{c}$.
               The units of $q_{c}$ are is $g/kg$.}
      \label{tab:qc0}
    \end{table}
    \begin{table}[H]
      \centering
      \begin{tabular}{ c| c c c c c c c c c}
        \hline
        \hline
                          & $st=1$  &  $st=2$ & $st=3$  && $st=1$ &  $st=2$ & $st=3$  &&\\
                          & \multicolumn{3}{c}{$\mathcal{P}_{5}$} && \multicolumn{3}{c}{$\mathcal{P}_{7}$} && Target\\
                          & \multicolumn{7}{c}{DBI} &&  \\
        \hline
          maximum $q_{c}$ &  0.20   & 0.20    & 0.31    && 0.30   & 0.30    & 0.28  && 0.28  \\ 
          total $q_{c}$   &  45.91  & 47.74   & 87.98   && 86.57  & 82.67   & 75.11 && 69.82\\
                          & \multicolumn{7}{c}{PPI}  &&  \\
        \hline
         maximum $q_{c}$ & 0.20   & 0.21   & 0.33   && 0.32  & 0.29   & 0.29   && 0.28 \\
         total $q_{c}$   & 47.87  & 50.09  & 97.60  && 92.54 & 81.44  & 78.85  && 69.82 \\ 
        \hline
        \hline
      \end{tabular}                                                                                  
      \caption{Maximum values of $q_{c}$ and the total amount of the cloud mixing ratio at $t = 5 h$ with $nz=600$ points and $\epsilon_{0}=\epsilon_{1}=10^{-5}$.
               The total amount of the cloud mixing ratio is calculated by estimating the integral $q_{c}$.
               The units of $q_{c}$ are is $g/kg$.}
      \label{tab:qc1}
    \end{table}

  \subsection{Example VII $f_7(x)$}
  \label{subsec:example5}
  This example uses an extended version of the 1D Runge function defined in Equation \ref{eq:f1} from Section \ref{subsec:example1} to a 2D function:  
  \begin{equation}\label{eq:f6}
    f_{7}(x, y) = \frac{1}{1+25(x^{2}+y^{2})}, \quad x, y \in [-1,1]
  \end{equation}
  For the DBI and PPI algorithm, $\epsilon_{0}=0.01$, $\epsilon_{1}=1$ and $st=2$.
  The results from Tables \ref{tab:runge2D_1} and \ref{tab:runge2Dp2} show that the DBI and PPI methods give smaller approximation errors compared to the PCHIP and SPS methods.
  In this case, the DBI and PPI methods use higher order polynomial interpolants for each interval.
  These higher order interpolants help improve the approximation compared to the PCHIP and  SPS.

  \begin{table}[H]
     \centering
     \scalebox{1.0}{
     \begin{tabular}{| c| c| c| c c| c c|}
          \hline
        $N_{i}$   & PCHIP    &SPS    & \multicolumn{2}{c|}{DBI} & \multicolumn{2}{c|}{PPI}     \\
          \hline              
                  &$L^{2}$-error & $L^{2}$-error & $L^{2}$-error & avg. deg. & $L^{2}$-error & avg. deg. \\
          \hline                   
                  &          &       & \multicolumn{4}{c|}{$\mathcal{P}_{1}$} \\
          \hline               
        $17$      &--        & --    & 1.60E-2       & 1      & 1.60E-2   & 1     \\
        $33$      &--        & --    & 4.42E-3       & 1      & 4.42E-3   & 1     \\
	$65$      &--        & --    & 1.12E-3       & 1      & 1.12E-3   & 1     \\
	$129$     &--        & --    & 2.82E-4       & 1      & 2.82E-4   & 1     \\
	$257$     &--        & --    & 7.06E-5       & 1      & 7.06E-5   & 1     \\
          \hline                   
                  &  $\mathcal{P}_{3}$ &   \multicolumn{5}{c|}{$\mathcal{P}_{4}$} \\
          \hline               
        $17$      & 5.01E-3  & 3.61E-3  & 5.16E-3   & 3.97   & 4.28E-3   & 4   \\
        $33$      & 1.23E-3  & 5.44E-4  & 3.51E-4   & 3.98   & 3.31E-4   & 4   \\
	$65$      & 2.33E-4  & 1.23E-4  & 2.55E-5   & 3.98   & 1.31E-5   & 4   \\
	$129$     & 4.27E-5  & 3.07E-5  & 1.20E-6   & 3.99   & 4.36E-7   & 4   \\
	$257$     & 7.72E-6  & 3.34E-6  & 4.96E-8   & 4      & 1.39E-8   & 4   \\
          \hline                                                                    
                  &          &        \multicolumn{5}{c|}{$\mathcal{P}_{8}$} \\
          \hline                          
        $17$      & --       & 8.55E-3  &  3.19E-3   & 7.75   & 1.84E-3   & 7.99  \\
        $33$      & --       & 1.99E-3  &  2.78E-4   & 7.82   & 7.86E-5   & 8  \\
        $65$      & --       & 4.97E-4  &  2.31E-5   & 7.90   & 5.14E-7   & 8  \\
        $129$     & --       & 1.25E-4  &  1.13E-6   & 7.95   & 1.49E-9   & 8  \\
        $257$     & --       & 3.16E-5  &  4.78E-8   & 7.98   & 3.25E-12  & 8  \\
          \hline                                                                    
                  &          &        \multicolumn{5}{c|}{$\mathcal{P}_{16}$} \\
          \hline                           
        $17$      & --       & 1.19E-2  &  3.49E-3   & 13.15   & 2.83E-3   & 14.26\\
        $33$      & --       & 3.10E-3  &  2.74E-4   & 15.43   & 2.68E-5   & 16\\
        $65$      & --       & 7.82E-4  &  2.30E-5   & 15.69   & 2.63E-8   & 16\\
        $129$     & --       & 1.96E-4  &  1.13E-6   & 15.84   & 1.77E-12   & 16\\
        $257$     & --       & 4.93E-5  &  4.76E-8   & 15.92   & 1.89E-15   & 16\\
        \hline
     \end{tabular}}                                                                                  
     \caption{$L^{2}-errors$ when approximating $f_{7}(x,y)$ with $N_{i}\times N_{i}$ points.
           $N_{i}$ represents the number of input points used in each dimension to build the approximation.
           $P_{j}$ represents the use of polynomials of degree $j$, with $j$ being the target degree.
           The fourth and sixth columns show the average polynomial degree used for the DBI and PPI methods, respectively.
           The mesh points are uniformly distributed on each dimension.}
     \label{tab:runge2D_1}
   \end{table}
   \begin{table}[H]
     \centering
     \scalebox{1.0}{
     \begin{tabular}{| c| c| c| c c| c c|}
          \hline
        $N_{i}$   & PCHIP   & SPS      & \multicolumn{2}{c|}{DBI} & \multicolumn{2}{c|}{PPI}     \\
          \hline             
                  &$L^{2}$-error & $L^{2}$-error & $L^{2}$-error & avg. deg. & $L^{2}$-error & avg. deg. \\
          \hline                   
                  &          &       & \multicolumn{4}{c|}{$\mathcal{P}_{1}$} \\
          \hline               
        $17$      &--        & --    & 1.60E-2       & 1      & 1.60E-2   & 1     \\
        $33$      &--        & --    & 4.42E-3       & 1      & 4.42E-3   & 1     \\
	$65$      &--        & --    & 1.12E-3       & 1      & 1.12E-3   & 1     \\
	$129$     &--        & --    & 2.82E-4       & 1      & 2.82E-4   & 1     \\
	$257$     &--        & --    & 7.11E-5       & 1      & 7.11E-5   & 1     \\
          \hline                    
                  &  $\mathcal{P}_{3}$  &     \multicolumn{5}{c|}{$\mathcal{P}_{4}$} \\
          \hline                
        $17$      & 6.31E-03 & 5.41E-3  & 5.54E-3   & 3.97   & 5.33E-3   & 4     \\
        $33$      & 4.15E-04 & 1.05E-3  & 4.59E-4   & 3.98   & 4.53E-4   & 4     \\
	$65$      & 1.07E-04 & 1.61E-4  & 2.67E-5   & 3.98   & 2.54E-5   & 4     \\
	$129$     & 2.46E-05 & 4.31E-5  & 7.97E-7   & 3.99   & 6.86E-7   & 4     \\
	$257$     & 5.20E-06 & 1.12E-5  & 2.79E-8   & 4.00   & 2.15E-8   & 4     \\
          \hline                                                                    
                  &         &           \multicolumn{5}{c|}{$\mathcal{P}_{8}$}\\
          \hline                          
        $17$      & --      & 5.89E-3  &  3.02E-3   & 7.75   & 2.81E-3    & 7.99    \\
        $33$      & --      & 1.72E-3  &  9.41E-5   & 7.82   & 9.34E-5    & 8    \\
        $65$      & --      & 5.44E-4  &  1.78E-6   & 7.90   & 1.51E-6    & 8    \\
        $129$     & --      & 1.90E-4  &  4.31E-8   & 7.95   & 4.53E-9    & 8    \\
        $257$     & --      & 4.91E-5  &  1.73E-9   & 7.98   & 1.87E-11   & 8    \\
          \hline                                                                       
                  &         &           \multicolumn{5}{c|}{$\mathcal{P}_{16}$}\\
          \hline                          
        $17$      & --      & 1.88E-2  &  6.31E-3   & 13.15   & 8.93E-3    & 14.26 \\
        $33$      & --      & 2.11E-3  &  2.93E-5   & 15.43   & 2.91E-5    & 16    \\
        $65$      & --      & 6.96E-4  &  8.12E-7   & 15.69   & 5.41E-8    & 16    \\
        $129$     & --      & 2.25E-4  &  2.27E-8   & 15.84   & 1.97E-11   & 16    \\
        $257$     & --      & 7.93E-5  &  9.21E-10  & 15.92   & 7.22E-15   & 16    \\
        \hline
     \end{tabular}}                                                                                  
     \caption{$L^{2}-errors$ when approximating $f_{7}(x,y)$ with $N_{i}\times N_{i}$ points.
           $N_{i}$ represents the number of input points used in each dimension to build the approximation.
           $P_{j}$ represents the use of polynomials of degree $j$, with $j$ being the target degree.
           The fourth and sixth columns show the average polynomial degree used for the DBI and PPI methods respectively.
           For each dimension, the interval $[-1,1]$ is divided into $(N_{i}-1) / j$ elements and $j+1$ LGL quadrature points are used in each element.}
     \label{tab:runge2Dp2}
   \end{table}

  \subsection{Example VIII $f_8(x)$}
  \label{subsec:example6}
  This example uses a 2D function that is used to study positive and monotonic splines \cite{CHAN2001135, 10.1007/11537908_20, lancaster1986curve}.
  The function is defined as follows:
  \begin{equation}\label{eq:T2}
    f_{8}(x, y) =
    \begin{cases}
      2(y-x) & \textrm{ if } 0 \leq y-x \leq 0.5 \\
      1      & \textrm{ if } y-x \geq 0.5 \\
      cos \bigg( 4\pi \sqrt{(x-1.5)^{2}+(y-0.5)^{2}}\bigg) & \textrm{ if }(x-1.5)^{2}+(y-0.5)^{2} \leq \frac{1}{16} \\
      0      & otherwise
    \end{cases}
  \end{equation}
  For the DBI and PPI algorithm, $\epsilon_{0}=0.01$, $\epsilon_{1}=1$ and $st=2$. 
  As in Example V, the function $f_{8}(x)$ is $C^{0}$-continuous and the underlying mesh used for the approximations does not capture the sharp corners.
  The $L^{2}$-errors from the DBI and PPI methods are dominated by the local errors of the intervals with $C^{0}$-continuity and low degree polynomial interpolants.
  Tables \ref{tab:T2_1} and \ref{tab:T2_2} show that the $L^{2}$-errors from the three methods have the same order,
with DBI and PPI having slightly smaller errors than the other approaches.
  In the cases where the underlying function is $C^{0}$, the results from DBI and PPI are comparable to the other approaches.
  Furthermore, the results from DBI and PPI can be improved by using a mesh that captures $C^{0}$-continuity, as is the case with the spectral element methods in NEPTUNE. 

  \begin{table}[H]
     \centering
     \begin{tabular}{| c| c| c| c c| c c|}
          \hline
        $N_{i}$   & PCHIP   & SPS      & \multicolumn{2}{c|}{DBI} & \multicolumn{2}{c|}{PPI}     \\
          \hline             
                  &$L^{2}$-error & $L^{2}$-error & $L^{2}$-error & avg. deg. & $L^{2}$-error & avg. deg. \\
          \hline             
                  &         &          & \multicolumn{4}{c|}{$\mathcal{P}_{1}$} \\
          \hline               
        $17$      & --      & --       & 2.70E-2     & 1   & 2.70E-2       & 1    \\
        $33$      & --      & --       & 9.51E-3     & 1   & 9.51E-3       & 1    \\
        $65$      & --      & --       & 3.40E-3     & 1   & 3.40E-3       & 1    \\
        $129$     & --      & --       & 1.20E-3     & 1   & 1.20E-3       & 1    \\
        $257$     & --      & --       & 4.30E-4     & 1   & 4.30E-4       & 1    \\
          \hline             
                  &  $\mathcal{P}_{3}$ &     \multicolumn{5}{c|}{$\mathcal{P}_{4}$} \\
          \hline               
        $17$      & 1.91E-2 & 1.87E-2  & 1.77E-2   & 2.04   & 1.73E-2   & 2.06    \\
        $33$      & 6.92E-3 & 6.11E-3  & 6.22E-3   & 1.93   & 6.21E-3   & 1.95    \\
        $65$      & 2.47E-3 & 2.69E-3  & 2.24E-3   & 1.89   & 2.24E-3   & 1.90    \\
        $129$     & 8.99E-4 & 7.71E-4  & 8.17E-4   & 1.88   & 8.16E-4   & 1.88    \\
        $257$     & 3.23E-4 & 2.77E-4  & 2.95E-4   & 1.87   & 2.94E-4   & 1.87    \\
          \hline                                                                 
                  &         &         \multicolumn{5}{c|}{$\mathcal{P}_{8}$} \\
          \hline                         
        $17$      & --      & 1.91E-2  & 1.73E-2   & 3.19   & 1.69E-2   & 3.31    \\
        $33$      & --      & 6.46E-3  & 6.20E-3   & 3.09   & 6.19E-3   & 3.16    \\
        $65$      & --      & 2.24E-3  & 2.21E-3   & 3.04   & 2.20E-3   & 3.09    \\
        $129$     & --      & 8.12E-4  & 7.98E-4   & 3.03   & 7.97E-4   & 3.04    \\
        $257$     & --      & 2.92E-4  & 2.87E-4   & 3.01   & 2.87E-4   & 3.02    \\
          \hline                                                                 
                  &         &         \multicolumn{5}{c|}{$\mathcal{P}_{16}$} \\
          \hline                         
        $17$      & --      & 2.19E-2  & 2.02E-2   & 4.58   & 2.34E-2   & 4.94    \\
        $33$      & --      & 7.57E-3  & 6.14E-3   & 5.13   & 6.17E-3   & 5.35    \\
        $65$      & --      & 2.68E-3  & 2.25E-3   & 5.25   & 2.26E-3   & 5.38    \\
        $129$     & --      & 9.63E-4  & 8.07E-4   & 5.29   & 8.08E-4   & 5.35    \\
        $257$     & --      & 3.45E-4  & 2.89E-4   & 5.29   & 2.89E-4   & 5.32    \\
        \hline
     \end{tabular}
     \caption{$L^{2}-errors$ when approximating $f_{8}(x,y)$ with $N_{i}\times N_{i}$ points.
           $N_{i}$ represents the number of input points used in each dimension to build the approximation.
           $P_{j}$ represents the use of polynomials of degree $j$, with $j$ being the target degree.
           The fourth and sixth columns show the average polynomial degree used for the DBI and PPI methods, respectively.
           The points are uniformly distributed in each dimension.}
     \label{tab:T2_1}
   \end{table}
   \begin{table}[H]
     \centering
     \begin{tabular}{| c| c| c| c c| c c|}
          \hline
        $N_{i}$   & PCHIP   &          & \multicolumn{2}{c|}{DBI} & \multicolumn{2}{c|}{PPI}     \\
          \hline             
                  &$L^{2}$-error & $L^{2}$-error & $L^{2}$-error & avg. deg. & $L^{2}$-error & avg. deg. \\
          \hline             
                  &         &          & \multicolumn{4}{c|}{$\mathcal{P}_{1}$} \\
          \hline                
        $17$      & --      & --       & 2.70E-2     & 1   & 2.70E-2       & 1    \\
        $33$      & --      & --       & 9.51E-3     & 1   & 9.51E-3       & 1    \\
        $65$      & --      & --       & 3.40E-3     & 1   & 3.40E-3       & 1    \\
        $129$     & --      & --       & 1.20E-3     & 1   & 1.20E-3       & 1    \\
        $257$     & --      & --       & 4.30E-4     & 1   & 4.30E-4       & 1    \\
          \hline             
                  & $\mathcal{P}_{3}$  &   \multicolumn{5}{c|}{$\mathcal{P}_{4}$} \\
          \hline                
        $17$      & 1.91E-2 & 2.55E-2  & 2.18E-2   & 2.04   & 2.17E-2   & 2.06   \\
        $33$      & 6.92E-3 & 5.76E-3  & 7.22E-3   & 1.93   & 7.18E-3   & 1.95   \\
        $65$      & 2.47E-3 & 2.11E-3  & 2.74E-3   & 1.89   & 2.72E-3   & 1.90   \\
        $129$     & 8.99E-4 & 8.08E-4  & 9.93E-4   & 1.88   & 9.87E-4   & 1.88   \\
        $257$     & 3.23E-4 & 2.92E-4  & 3.61E-4   & 1.87   & 3.59E-4   & 1.87   \\
          \hline                                                                 
                  &         &  \multicolumn{5}{c|}{$\mathcal{P}_{8}$} \\
          \hline                         
        $17$      & --      & 4.06E-2  & 3.42E-2   & 3.19   & 3.19E-2   & 3.31   \\
        $33$      & --      & 9.69E-3  & 8.68E-3   & 3.09   & 8.67E-3   & 3.16   \\
        $65$      & --      & 2.46E-3  & 2.76E-3   & 3.04   & 2.82E-3   & 3.09   \\
        $129$     & --      & 9.83E-4  & 1.09E-3   & 3.03   & 1.12E-3   & 3.04   \\
        $257$     & --      & 3.63E-4  & 3.85E-4   & 3.01   & 3.99E-4   & 3.02   \\
          \hline                                                                 
                  &         &  \multicolumn{5}{c|}{$\mathcal{P}_{16}$} \\
          \hline                         
        $17$      & --      & 4.36E-2  & 3.35E-2   & 4.58   & 2.81E-2   & 4.94   \\
        $33$      & --      & 1.53E-2  & 1.06E-2   & 5.13   & 1.06E-2   & 5.35   \\
        $65$      & --      & 4.31E-3  & 3.42E-3   & 5.25   & 3.46E-3   & 5.38   \\
        $129$     & --      & 1.13E-3  & 9.84E-4   & 5.29   & 1.29E-3   & 5.35   \\
        $257$     & --      & 4.38E-4  & 3.95E-4   & 5.29   & 5.07E-4   & 5.32   \\
        \hline
     \end{tabular}
     \caption{$L^{2}-errors$ when approximating $f_{8}(x,y)$ with $N_{i}\times N_{i}$ points.
           $N_{i}$ represents the number of input points used in each dimension to build the approximation.
           $P_{j}$ represents the use of polynomials of degree $j$, with $j$ being the target degree.
           The fourth and sixth columns show the average polynomial degree used for the DBI and PPI methods respectively.
           For each dimension, the interval $[-1,1]$ is divided into $(N_{i}-1) / j$ elements and $j+1$ LGL quadrature points are used in each element.}
     \label{tab:T2_2}
   \end{table}

  \subsection{Example IX  $f_9(x)$}
  \label{subsec:example8}
  This example is used herein to study shape-preserving (monotonicity and convexity) splines  \cite{10.2307/2157975}. 
  \begin{equation}\label{eq:T1}
    f_{9}(x,y) = max \bigg(0, sin(\pi x)sin(\pi y) \bigg) \quad x,y \in [-1,1]
  \end{equation}
  For the DBI and PPI algorithm, $\epsilon_{0}=0.01$, $\epsilon_{1}=1$ and $st=2$. 
  The function $f_{9}(x, y)$ is a $C^{0}$-continuous function.
  Tables \ref{tab:T1_1} and \ref{tab:T2_2} show $L^{2}$-errors when approximating $f_{9}(x,y)$ with the PCHIP, SPS, DBI, and PPI methods.
  The underlying mesh is such that the $C^{0}$-continuities are at the elements boundaries except for $\mathcal{P}_{16}$ and $N=17$.
  The PCHIP and SPS methods struggle to capture the $C^{0}$-continuities because both methods enforce $C^{1}$-continuity.  
  The $L^{2}$-error from DBI is dominated by the local error from the intervals with low-degree interpolants and so 
  as the average polynomial degree increases the $L^{2}$-errors do not improve.
  The $L^{2}$-error for $\mathcal{P}_{16}$ and $N=17$ is larger compared to the other cases when the PPI method is used.
  For $\mathcal{P}_{16}$ and $N=17$, there  is no mesh point at the points of $C^{0}$-continuity and so the $L^{2}$-error is dominated by the local error from those intervals where low degree interpolants are used.
  Overall, the results from Tables \ref{tab:T1_1} and \ref{tab:T2_2} demonstrate that the DBI and PPI methods lead to smaller approximation errors than the PCHIP and  SPS methods.
 \begin{table}[H]
    \centering
    \scalebox{1.0}{
    \begin{tabular}{| c| c| c| c c| c c|}
          \hline
        $N_{i}$   & PCHIP    & SPS      & \multicolumn{2}{c|}{DBI} & \multicolumn{2}{c|}{PPI}     \\
          \hline              
                  &$L^{2}$-error & $L^{2}$-error & $L^{2}$-error & avg. deg. & $L^{2}$-error & avg. deg. \\
          \hline               
                  &          &          & \multicolumn{4}{c|}{$\mathcal{P}_{1}$} \\
          \hline                         
        $17$      & --      & --        & 5.80E-2  & 1     & 5.80E-2  & 1    \\
        $33$      & --      & --        & 1.88E-2  & 1     & 1.88E-2  & 1    \\
        $65$      & --      & --        & 6.27E-3  & 1     & 6.27E-3  & 1    \\
        $129$     & --      & --        & 2.17E-3  & 1     & 2.17E-3  & 1    \\
        $257$     & --      & --        & 7.87E-4  & 1     & 7.87E-4  & 1    \\
          \hline               
                  & $\mathcal{P}_{3}$   &         \multicolumn{5}{c|}{$\mathcal{P}_{4}$} \\
          \hline                         
        $17$      & 1.91E-2  & 1.80E-2  & 1.20E-2   & 2.56   & 1.20E-2   & 2.66  \\
        $33$      & 6.77E-3  & 6.21E-3  & 4.25E-3   & 2.53   & 4.25E-3   & 2.58  \\
        $65$      & 2.39E-3  & 2.18E-3  & 1.50E-3   & 2.52   & 1.50E-3   & 2.54  \\
        $129$     & 8.47E-4  & 7.70E-4  & 5.30E-4   & 2.51   & 5.30E-4   & 2.52  \\
        $257$     & 3.01E-4  & 2.74E-4  & 1.88E-4   & 2.50   & 1.88E-4   & 2.51  \\
          \hline                                                                 
                  &          &         \multicolumn{5}{c|}{$\mathcal{P}_{8}$} \\
          \hline                          
        $17$      & --       & 1.47E-2  & 1.27E-2   & 4.22   & 1.27E-2   & 4.72  \\
        $33$      & --       & 4.57E-3  & 4.48E-3   & 4.37   & 4.48E-3   & 4.61  \\
        $65$      & --       & 1.51E-4  & 1.58E-3   & 4.44   & 1.58E-3   & 4.56  \\
        $129$     & --       & 5.16E-4  & 5.60E-4   & 4.47   & 5.60E-4   & 4.53  \\
        $257$     & --       & 1.84E-4  & 1.98E-4   & 4.48   & 1.98E-4   & 4.51  \\
          \hline                                                                 
                  &          &         \multicolumn{5}{c|}{$\mathcal{P}_{16}$} \\
          \hline                          
        $17$      & --       & 1.50E-2  & 2.22E-2   & 5.92   & 2.74E-2   & 7.11  \\
        $33$      & --       & 4.03E-3  & 4.79E-3   & 8.05   & 4.79E-3   & 8.67  \\
        $65$      & --       & 1.15E-3  & 1.69E-3   & 8.28   & 1.69E-3   & 8.58  \\
        $129$     & --       & 3.50E-4  & 5.98E-4   & 8.39   & 5.98E-4   & 8.54  \\
        $257$     & --       & 1.17E-4  & 2.12E-4   & 8.44   & 2.12E-4   & 8.52  \\
          \hline              
     \end{tabular}}                                                                                  
     \caption{$L^{2}-errors$ when approximating $f_{9}(x,y)$ with $N_{i}\times N_{i}$ points.
           $N_{i}$ represents the number of input points used in each dimension to build the approximation.
           $P_{j}$ represents the use of polynomials of degree $j$, with $j$ being the target degree.
           The fourth and sixth columns show the average polynomial degree used for the DBI and PPI methods respectively.
           The points are uniformly distributed on each dimension.}
     \label{tab:T1_1}
  \end{table}
  \begin{table}[H]
     \centering
     \scalebox{1.0}{
     \begin{tabular}{| c| c| c| c c| c c|}
          \hline
        $N_{i}$   & PCHIP    & SPS      & \multicolumn{2}{c|}{DBI} & \multicolumn{2}{c|}{PPI}     \\
          \hline              
                  &$L^{2}$-error & $L^{2}$-error & $L^{2}$-error & avg. deg. & $L^{2}$-error & avg. deg. \\
          \hline              
                  &          &          & \multicolumn{4}{c|}{$\mathcal{P}_{1}$} \\
          \hline                 
        $17$      & --      & --        & 5.80E-2  & 1     & 5.80E-2  & 1    \\
        $33$      & --      & --        & 1.88E-2  & 1     & 1.88E-2  & 1    \\
        $65$      & --      & --        & 6.27E-3  & 1     & 6.27E-3  & 1    \\
        $129$     & --      & --        & 2.17E-3  & 1     & 2.17E-3  & 1    \\
        $257$     & --      & --        & 7.87E-4  & 1     & 7.87E-4  & 1    \\
          \hline              
                  & $\mathcal{P}_{3}$    &        \multicolumn{5}{c|}{$\mathcal{P}_{4}$} \\
          \hline                 
        $17$      & 3.56E-02& 1.05E-2  & 1.85E-2  & 2.52  & 1.46E-4  & 3.91     \\
        $33$      & 1.79E-03& 3.61E-3  & 4.74E-3  & 2.61  & 4.77E-6  & 3.95     \\
        $65$      & 2.53E-04& 1.26E-3  & 1.19E-3  & 2.62  & 1.43E-7  & 3.98     \\
        $129$     & 9.94E-05& 4.44E-4  & 2.98E-4  & 2.62  & 9.16E-9  & 3.99     \\
        $257$     & 3.26E-05& 1.61E-4  & 7.45E-5  & 2.62  & 5.52E-9  & 3.99     \\
          \hline              
                  &          &        \multicolumn{5}{c|}{$\mathcal{P}_{8}$} \\
          \hline                 
        $17$      & --       & 1.61E-2  & 1.85E-2  & 4.04  & 6.01E-8  & 7.36     \\
        $33$      & --       & 3.34E-3  & 4.74E-3  & 4.57  & 1.37E-8  & 7.89     \\
        $65$      & --       & 8.74E-4  & 1.19E-3  & 4.57  & 9.43E-9  & 7.95     \\
        $129$     & --       & 2.40E-4  & 2.98E-4  & 4.56  & 6.15E-9  & 7.97     \\
        $257$     & --       & 7.32E-5  & 7.45E-5  & 4.56  & 3.57E-9  & 7.99     \\
          \hline              
                  &          &        \multicolumn{5}{c|}{$\mathcal{P}_{16}$} \\
          \hline                 
        $17$      & --       & 2.72E-2  & 1.85E-2  & 4.29  & 2.35E-3  & 8.75     \\
        $33$      & --       & 6.65E-3  & 4.74E-3  & 7.93  & 9.88E-9  & 15.31    \\
        $65$      & --       & 1.33E-3  & 1.19E-3  & 8.57  & 6.44E-9  & 15.88    \\
        $129$     & --       & 3.35E-4  & 2.98E-4  & 8.54  & 3.78E-9  & 15.94    \\
        $257$     & --       & 8.41E-5  & 7.45E-5  & 8.55  & 1.75E-9  & 15.97    \\
          \hline              
     \end{tabular}}                                                                                  
     \caption{$L^{2}-errors$ when approximating $f_{9}(x,y)$ with $N_{i}\times N_{i}$ points.
           $N_{i}$ represents the number of input points used in each dimension to build the approximation.
           $P_{j}$ represents the use of polynomials of degree $j$, with $j$ being the target degree.
           The fourth and sixth columns show the average polynomial degree used for the DBI and PPI methods respectively.
           For each dimension, the interval $[-1,1]$ is divided into $(N_{i}-1) / j$ elements 
           and $j+1$ LGL quadrature points are used in each element.}
     \label{tab:T1_2}
  \end{table}

  \subsection{Example X $f_{10}(x)$ }
  \label{subsec:example9}
  This example uses a 2D extension of the 1D approximation of the Heaviside function $f_{2}(x)$ defined in Equation \ref{eq:f2} which
  is defined as follows:
  \begin{equation}\label{eq:Heaviside2}
    f_{10}(x,y) = \frac{1}{1+e^{-\sqrt{2}k(x+y)}}, \quad x, y \in [-0.2, 0.2]
  \end{equation}
  For the DBI and PPI algorithm, $\epsilon_{0}=0.01$, $\epsilon_{1}=1$ and $st=2$. 
  The function $f_{10}(x,y)$ is challenging because of the large gradient at $y=-x$.
  Tables \ref{tab:Heaviside2D_1} and \ref{tab:Heaviside2D_2} show $L^{2}$-errors when approximating $f_{10}(x)$ using PCHIP, SPS, DBI, and PPI.
  As the average polynomial degree increases the accuracy of the DBI and PPI methods improves.
  In this case, the $L^{2}$-error is dominated by the local error of the region with the steep gradient.
  The errors for the DBI and PPI methods are similar because the stencils used for both methods are the same in the region with the large gradient. 
  Overall, the results from the Tables \ref{tab:Heaviside2D_1} and \ref{tab:Heaviside2D_2} show that the DBI and PPI methods lead to smaller $L^{2}$-errors compared to the other methods.
  \begin{table}[H]
     \centering
     \scalebox{1.0}{
     \begin{tabular}{| c| c| c| c c| c c|}
          \hline
        $N_{i}$   & PCHIP    & SPS      & \multicolumn{2}{c|}{DBI} & \multicolumn{2}{c|}{PPI}     \\
          \hline              
                  &$L^{2}$-error & $L^{2}$-error & $L^{2}$-error & avg. deg. & $L^{2}$-error & avg. deg. \\
          \hline              
                  &          &          & \multicolumn{4}{c|}{$\mathcal{P}_{1}$} \\
          \hline                 
        $17$      &  --      & --       & 1.50E-2 & 1      & 1.50E-2  & 1    \\
        $33$      &  --      & --       & 4.57E-3 & 1      & 4.57E-3  & 1    \\
        $65$      &  --      & --       & 1.26E-3 & 1      & 1.26E-3  & 1    \\
        $129$     &  --      & --       & 3.23E-4 & 1      & 3.23E-4  & 1    \\
        $257$     &  --      & --       & 8.15E-5 & 1      & 8.15E-5  & 1    \\
          \hline              
                  & $\mathcal{P}_{3}$  &   \multicolumn{5}{c|}{$\mathcal{P}_{4}$} \\
          \hline                 
        $17$      & 8.07E-3  & 1.99E-3  & 9.45E-3   & 2.75   & 9.44E-3   & 3.23   \\
        $33$      & 1.26E-3  & 2.43E-4  & 1.33E-3   & 3.53   & 1.29E-3   & 3.65   \\
        $65$      & 1.44E-4  & 4.92E-5  & 9.29E-5   & 3.82   & 9.29E-5   & 3.82   \\
        $129$     & 1.63E-5  & 1.20E-5  & 3.67E-6   & 3.81   & 3.67E-6   & 3.81   \\
        $257$     & 1.94E-6  & 3.05E-6  & 1.21E-7   & 3.79   & 1.21E-7   & 3.79   \\
          \hline                                                                  
                  &          &     \multicolumn{5}{c|}{$\mathcal{P}_{8}$} \\
          \hline                          
        $17$      &  --      & 1.80E-1  & 8.05E-3   & 3.15   & 8.67E-3   & 5      \\
        $33$      &  --      & 1.22E-1  & 1.03E-3   & 5.78   & 9.05E-4   & 6.55   \\
        $65$      &  --      & 8.48E-2  & 4.83E-5   & 7.54   & 4.99E-5   & 7.57   \\
        $129$     &  --      & 1.04E-1  & 2.57E-7   & 7.53   & 2.57E-7   & 7.55   \\
        $257$     &  --      & 8.02E-2  & 5.27E-10  & 7.49   & 5.27E-10  & 7.52   \\
          \hline                                                                  
                  &          &     \multicolumn{5}{c|}{$\mathcal{P}_{16}$} \\
          \hline                          
        $17$      &  --      & 4.72E-3  & 7.39E-3   & 3.57   & 1.86E-2   & 7.45   \\
        $33$      &  --      & 1.22E-3  & 1.02E-3   & 6.71   & 2.33E-3   & 10.23  \\
        $65$      &  --      & 3.11E-4  & 2.12E-4   & 14.76  & 2.41E-4   & 14.94  \\
        $129$     &  --      & 7.80E-5  & 1.03E-6   & 14.93  & 1.03E-6   & 15.05  \\
        $257$     &  --      & 1.95E-5  & 4.41E-11  & 14.83  & 4.41E-11  & 14.97  \\
          \hline              
     \end{tabular}}                                                                                  
     \caption{$L^{2}-errors$ when approximating $f_{10}(x,y)$ with $N_{i}\times N_{i}$ points.
           $N_{i}$ represents the number of input points used in each dimension to build the approximation.
           $P_{j}$ represents the use of polynomials of degree $j$, with $j$ being the target degree.
           The fourth and sixth columns show the average polynomial degree used for the DBI and PPI methods, respectively.
           The points are uniformly distributed on each dimension.}
     \label{tab:Heaviside2D_1}
  \end{table}
  \begin{table}[H]
     \centering
     \scalebox{1.0}{
     \begin{tabular}{| c| c| c| c c| c c|}
          \hline
        $N_{i}$   & PCHIP    & SPS      & \multicolumn{2}{c|}{DBI} & \multicolumn{2}{c|}{PPI}     \\
          \hline              
                  &$L^{2}$-error & $L^{2}$-error & $L^{2}$-error & avg. deg. & $L^{2}$-error & avg. deg. \\
          \hline              
                  &          &          & \multicolumn{4}{c|}{$\mathcal{P}_{1}$} \\
          \hline                 
        $17$      &  --      & --       & 1.50E-2 & 1      & 1.50E-2  & 1    \\
        $33$      &  --      & --       & 4.57E-3 & 1      & 4.57E-3  & 1    \\
        $65$      &  --      & --       & 1.26E-3 & 1      & 1.26E-3  & 1    \\
        $129$     &  --      & --       & 3.23E-4 & 1      & 3.23E-4  & 1    \\
        $257$     &  --      & --       & 8.15E-5 & 1      & 8.15E-5  & 1    \\
          \hline              
                  & $\mathcal{P}_{3}$ &        \multicolumn{5}{c|}{$\mathcal{P}_{4}$} \\
          \hline                 
        $17$      &  1.32E-2 & 2.79E-3 & 1.12E-2   & 2.75   & 1.11E-2   & 3.23  \\
        $33$      &  2.75E-3 & 3.49E-4 & 1.73E-3   & 3.53   & 1.68E-3   & 3.65  \\
        $65$      &  3.57E-4 & 6.87E-5 & 1.28E-4   & 3.82   & 1.28E-4   & 3.82  \\
        $129$     &  4.09E-5 & 1.64E-5 & 5.47E-6   & 3.81   & 5.47E-6   & 3.81  \\
        $257$     &  5.04E-6 & 4.12E-6 & 1.77E-7   & 3.79   & 1.77E-7   & 3.79  \\
          \hline                                                                 
                  &          &       \multicolumn{5}{c|}{$\mathcal{P}_{8}$} \\
          \hline                          
        $17$      &  --      & 5.82E-3 & 1.22E-2   & 3.15   & 1.20E-2   & 5.00  \\
        $33$      &  --      & 1.26E-3 & 1.82E-3   & 5.78   & 1.67E-3   & 6.55  \\
        $65$      &  --      & 3.00E-4 & 4.98E-5   & 7.54   & 4.98E-5   & 7.57  \\
        $129$     &  --      & 7.58E-5 & 4.03E-7   & 7.53   & 4.03E-7   & 7.55  \\
        $257$     &  --      & 1.90E-5 & 1.21E-9   & 7.49   & 1.21E-9   & 7.52  \\
          \hline                                                                 
                  &          &       \multicolumn{5}{c|}{$\mathcal{P}_{16}$} \\
          \hline                          
        $17$      &  --      & 8.05E-3 & 1.34E-2   & 3.57   & 1.31E-2   & 7.45  \\
        $33$      &  --      & 2.06E-3 & 2.04E-3   & 6.71   & 1.92E-3   & 10.23 \\
        $65$      &  --      & 5.14E-4 & 3.81E-5   & 14.76  & 3.84E-5   & 14.94 \\
        $129$     &  --      & 1.26E-4 & 6.20E-8   & 14.93  & 6.20E-8   & 15.05 \\
        $257$     &  --      & 3.18E-5 & 6.20E-12  & 14.83  & 6.20E-12  & 14.97 \\
          \hline              
     \end{tabular}}                                                                                  
     \caption{$L^{2}-errors$ when approximating $f_{10}(x,y)$ with $N_{i}\times N_{i}$ points.
           $N_{i}$ represents the number of input points used in each dimension to build the approximation.
           $P_{j}$ represents the use of polynomials of degree $j$, with $j$ being the target degree.
           The fourth and sixth columns show the average polynomial degree used for the DBI and PPI methods respectively.
           For each dimension, the interval $[-0.2,0.2]$ is divided into $(N_{i}-1) / j$ elements and $j+1$ LGL quadrature points are used in each element.}
     \label{tab:Heaviside2D_2}
  \end{table}

\section{Discussion and Conclusion}
In this report, a representative sample of existing methods is compared against our new approaches on
a number of different test functions, including smooth, $C^{0}$, discontinuous, and steep functions.
The comparison undertaken here focuses on how accurately the different methods are able to represent this underlying set of test functions.
Overall, the DBI and PPI methods perform well and are suited to the $C^{0}$ continuity of the spectral element methods in NEPTUNE.
The experiments show that the DBI and PPI methods are suitable approaches for interpolating
smooth functions and $C^0$ continuous functions while enforcing positivity.
In detail the summary is that:
\begin{itemize}
\item
The results in Section 3 Examples I, II, and III show that the improved DBI and new PPI approaches preserve positivity exactly as the proofs in \cite{ouermi2022eno} indicate;
\item
The results in Section 4 and Sections 5.1, 5.2, and 5.3 show that the DBI and PPI approaches give much higher levels of accuracy than the 
DBI method by allowing the solution to be outside the local bounds while remaining positive.
The PPI method also appears to give better results than the SPS method in line with  
the studies in \cite{Costantini1990} and \cite{10.1145/264029.264050} which demonstrate that the SPS method does not achieve high-order accuracy; and
\item 
In Examples V and VI, better approximations are obtained using $st=1$ and $st=3$, respectively, for the PPI and DBI methods.
These results show that in some cases the appropriate choice of the parameter $st$ can further improve the approximation in addition to preserving data-boundedness or positivity. 
Additional studies evaluating the different choices of $\epsilon_{0}$, $\epsilon_{1}$, and $st$ are presented in ~\cite{tajo2022PPIsoftware}.
\item
In the cases when steep gradients or discontinuities force the use of low-order approximations, the DBI and PPI methods compete against the well-known cubic spline method PCHIP and the higher order MQS and the SPS spline methods.
\end{itemize}
Overall, it would seem that when it is possible to use higher-order polynomial approximations the PPI method appears to give levels of accuracy 
that compete with standard unmodified high-order spline methods while at the same time preserving positivity.
\section*{Acknowledgements}
This work has been supported by the US Naval Research Laboratory (559000669), 
the National Science Foundation (1521748), and the Intel Graphics and Visualization Institute at the University of Utah's Scientific Computing and Imaging (SCI) 
Institute (29715).
The authors would like to thank  Dr. Alex Reinecke of the Naval Research Laboratory for his constant support and help.
\bibliographystyle{unsrt}
\bibliography{references2}
\end{document}